\documentclass[12pt]{article}
\usepackage{mathdots}
\usepackage{amsfonts}
\usepackage{amssymb}
\usepackage{latexsym}
\usepackage{graphicx}
\usepackage{pgfplots}
\usepackage{todonotes}

\usepackage{color}
\usepackage{fullpage}
\usepackage{verbatim}
\usepackage[section]{algorithm}
 \usepackage{algpseudocode}
\usepackage{algorithmicx}
\usepackage{hyperref}
\usepackage{amsopn,url}
\usepackage{amsmath}
\usepackage{amstext,amsfonts,bm,amssymb}
\usepackage{cleveref}
\newtheorem{theorem}{Theorem}[section]
\newtheorem{lemma}[theorem]{Lemma}

\newtheorem{corollary}[theorem]{Corollary}

 \numberwithin{equation}{section}
\algrenewcommand\algorithmicrequire{\textbf{Input:}}
\algrenewcommand\algorithmicensure{\textbf{Output:}}

\DeclareMathOperator*{\argmin}{argmin}

\newcommand*{\Scale}[2][4]{\scalebox{#1}{\ensuremath{#2}}}%
       \newcommand{\rowsize}{.7cm}
       \newcommand{\rrowsize}{.5cm}
       \newcommand{\colsize}{2.2cm}

       \newcommand{\smallsize}{.6cm}
       \newcommand{\bigsize}{1.2cm}

\def\Pi{\prod}

\newcommand{\uni}[1]{{\left\vert\kern-0.25ex\left\vert\kern-0.25ex\left\vert #1  \right\vert\kern-0.25ex\right\vert\kern-0.25ex\right\vert}} 
\newcommand{\unii}[1]{{\vert\kern-0.25ex\vert\kern-0.25ex\vert #1  \vert\kern-0.25ex\vert\kern-0.25ex\vert}} 
\newcommand{\Uniinv}[1]{{\big\vert\kern-0.25ex\big\vert\kern-0.25ex\big\vert #1  \big\vert\kern-0.25ex\big\vert\kern-0.25ex\big\vert}} 
\newcommand{\uniinv}[1]{{\left\vert\kern-0.25ex\left\vert\kern-0.25ex\left\vert #1  \right\vert\kern-0.25ex\right\vert\kern-0.25ex\right\vert}} 

\newcommand{\rr}[1]{\textcolor{red}{#1}}

\newcommand{\bb}[1]{\textcolor{blue}{#1}}

\newcommand{\ignore}[1]{}

\title{
SubApSnap: Solving parameter-dependent linear systems
with a snapshot
and subsampling
}
\author{
Eleanor Jones and Yuji Nakatsukasa\thanks{Mathematical Institute, University of Oxford. {\tt 
eleanor.jones@chch.ox.ac.uk}, {\tt yuji.nakatsukasa@maths.ox.ac.uk}\ This work is supported by EPSRC grants EP/Y010086/1 and EP/Y030990/1.
}
}

\begin{document}
\maketitle


\begin{abstract}
A growing number of problems in computational mathematics can be reduced to the solution of many linear systems that are related, often  depending smoothly or slowly on a parameter $p$, that is, $A(p)x(p)=b(p)$. 
We introduce an efficient algorithm for solving such parameter-dependent linear systems for many 
values of $p$. The algorithm, which we call SubApSnap (for \emph{Sub}sampled $A(p)$ times \emph{Snap}shot), is based on combining ideas from model order reduction and randomised linear algebra: namely, taking a snapshot matrix, and solving the resulting tall-skinny least-squares problems using a subsampling-based dimension-reduction approach. 
We show that SubApSnap is a strict generalisation of the popular DEIM algorithm in nonlinear model order reduction. 
SubApSnap is a sublinear-time algorithm, and 
once the snapshot and subsampling are determined, it solves $A(p_*)x(p_*)=b(p_*)$ for a new value of $p_*$ at a dramatically improved speed:
it does not even need to read the whole matrix $A(p_*)$ to solve the linear system for a new value of $p_*$. 
We prove under natural assumptions that, given a good subsampling and snapshot, SubApSnap yields solutions with small residual for all parameter values of interest. 
We illustrate the efficiency and performance of the algorithm with problems arising in PDEs, model reduction, and kernel ridge regression, where SubApSnap achieves speedups of many orders of magnitude over a standard solution; for example over $20,000\times$ for a $10^7\times 10^7$ problem, while providing good accuracy. 

\end{abstract}

\section{Introduction}
Solving a large-scale linear system $Ax=b$ 
for $A\in\mathbb{R}^{n\times n}, b\in\mathbb{R}^{n}$
continues to be among the most prevalent problems in computational mathematics. Many problems naturally lead to $Ax=b$, and even for nonlinear problems, a standard approach is to iteratively find a linear approximation and solve the linear problem, as in e.g. Newton's method. Classical methods and origins of $Ax=b$ (e.g. PDEs) can be found in Saad's textbook~\cite{saad2003iterative}. 

With the rise of data science and machine learning, linear systems 
have found new uses and applications. 
These include Gaussian processes~\cite{rasmussen2006gaussian}, ridge regression~\cite{friedman2001elements}, 
kernel methods~\cite{shawe2004kernel}, 
and model-order reduction~\cite{chaturantabut2010nonlinear}.
Increasingly commonly, it is not just a single linear system but a sequence or collection of them that are of interest, often depending on a parameter $p\in\Omega\subset \mathbb{R}^d$, 
that is, 
\begin{equation}\label{eq:problem}
A(p)x(p)=b(p),    
\end{equation} 
where the entries of $A(p),b(p)$ usually depend smoothly on $p\in \Omega$, e.g. analytically. 
For a fixed value of $p$,~\eqref{eq:problem} is simply an $n\times n$ linear system. 
Such problems arise for example in (hyper)parameter tuning~\cite{rasmussen2006gaussian}, model order reduction (e.g. where the transfer function $c^T(pE-A)^{-1}b$ needs to be evaluated at many values of $p$)~\cite{antoulas2020interpolatory,benner2015survey}, and parameter-dependent PDEs~\cite{barrault2004empirical,gunzburger2007reduced,hesthaven2016certified}. 
This paper describes an efficient numerical algorithm, SubApSnap, to find an approximate solution for $A(p_*)x(p_*)=b(p_*)$ for new/unseen value(s) of $p_*\in\Omega$. 
The set $\Omega$ is assumed to be bounded, and can lie in a high-dimensional space $d\geq 1$, i.e., $p$ may be a vector. 
More generally, SubApSnap can be used 
in situations where an explicit parameter is absent but a number of related linear systems need to be solved. 
For example, Newton's method for minimising $f(x)$ 
for a function $f:\mathbb{R}^{n}\rightarrow \mathbb{R}$
can be seen as an example where $A(p)$ is the Hessian of $f$ at $p=x\in\mathbb{R}^n$, so $d=n$. 
Nonetheless, it may be helpful to consider the simple 1-d case $\Omega=[-1,1]$, so $-1\leq p\leq 1$ on first reading. 



SubApSnap is based on two simple key ideas: \emph{Snapshot} and \emph{Subsampling}. 
We first form a subspace spanned by a snapshot, i.e., solutions for several different values of $p$. 
For a new parameter value $p_*$, 
finding a solution in that subspace leads to a least-squares problem; we then solve it via a carefully chosen subsampling strategy. 
Fortunately, for both snapshot and subsampling there is now a significant body of literature on their analysis and algorithms, namely in model order reduction for snapshots, and in (classical and randomized) linear algebra for subsampling. We discuss these further in Sections~\ref{sec:snap} and \ref{sec:subsample}. 

SubApSnap is thus able to find a solution for an $n\times n$ system $A(p_*)x(p_*)=b(p_*)$ in $O(nr^2)$ operations, where $r$ is the number of snapshot points taken; usually $r\ll n$, and often $r=O(1)$, so 
SubApSnap is a sublinear-time algorithm, i.e, less than $O(n^2)$ for dense problems. 
The reduction from the standard $O(n^3)$ in complexity is enormous. 
Further, SubApSnap can be even faster: for example when $A(p_*)$ is sparse with $O(1)$ nonzeros per row, its complexity becomes $O(r^3)$, so essentially $O(1)$. 
Crucially, SubApSnap does not even need to read the whole matrix $A(p_*)$, only a $O(r)$ subset of the rows. This can be particularly attractive when building the matrix $A(p_*)$ requires significant work, as is the case e.g. in kernel methods and some parameter-dependent PDEs. 
We demonstrate the accuracy and speed of SubApSnap through numerical experiments, where it is seen to achieve up to $30,000\times$ speedup over solving each linear system from scratch, while giving good accuracy. 

\section{Existing algorithms}
We will discuss methods that are closely related to  SubApSnap~\cite{guido2024subspace,kressner2011low,markovinovic2006accelerating}
in Section~\ref{sec:related}. 
Approaches that make use of previous solutions include Krylov subspace recycling~\cite{recyclekilmer,sturler06}, and rational Krylov methods~\cite{simoncini2010extended}. 
A Krylov method can also be truncated to take advantage of low-rank structures in the solution~\cite{kressner2011low}.
Recent studies~\cite{correnty2025chebyshev,correnty2023preconditioned,jarlebring2022infinite} approximate $A(p)$ and/or the solution using Taylor or Chebyshev expansions and/or grid sampling. 
These methods require that the system~\eqref{eq:problem} and the solution $x(p)$ depend smoothly (e.g. analytically) on $p$. 
Some of these are related to taking snapshots, but none employs subsampling, a crucial component of SubApSnap.

Snapshot is a popular approach in model order reduction~\cite{berkooz1993proper,chatterjee2000introduction}, and has been suggested previously for $A(p_*)x(p_*)=b(p_*)$~
\cite{markovinovic2006accelerating,guido2024subspace}.
Subsampling is a technique that has gained interest partly due to the advances of randomised algorithms in numerical linear algebra, and is a central topic in the closely related fields of active learning~\cite{roy2001toward} and experimental design~\cite{pukelsheim2006optimal}. 
Subsampling has also been used extensively in model reduction, namely in DEIM~\cite{chaturantabut2010nonlinear}, which also employs a snapshot, and we detail the similarities and differences between DEIM and SubApSnap in Section~\ref{sec:subapsnapvsDEIM}, where we show that SubApSnap can be seen as a  strict generalisation of DEIM.

It appears that methods that combine a snapshot and (randomised)SubApSnap is a  sketching have not been explored previously for parameter-dependent linear systems~\eqref{eq:problem}. In particular, most existing methods require reading 
an $n\times n$ matrix at each parameter value, so they are not sublinear time. Exceptions include those that explicitly find a representation of $x(p)$ with respect to $p$ using interpolation, e.g.~\cite{correnty2025chebyshev}, but such methods require strong assumptions on $x(p)$, such as smoothness, and struggle for example when $A(p)$ is (near-)singular for some $p\in\Omega$. 

{\it Notation}. 
For an $n\times r$ ($n\geq r$) matrix, its singular values of a matrix $A$ are $\sigma_i(A)$ for $1\leq i\leq r$, and $\sigma_{\max}(A)=\sigma_1(A),\sigma_{\min}(A)=\sigma_r(A)$ are the largest and smallest singular values.
$\|\cdot\|_2$ denotes the 2-norm (operator norm) of a matrix, or the Euclidean norm of a vector. 
For simplicity our discussion assumes $A(p)$ is a real-valued matrix, although everything extends to the complex case. 
\section{SubApSnap: Snapshot and subsampling for $A(p)x(p)=b(p)$}
The idea of snapshot applied to $A(p)x(p)=b(p)$ is as follows. 
In the initial, \emph{offline} phase, we solve $A(p_i)x(p_i)=b(p_i)$ (either to full precision or approximately) for prescribed values of $\{p_i\}_{i=1}^r$ with $p_i\in\Omega$, and form the snapshot matrix $X = [x(p_1),x(p_2),\ldots,x(p_r)]\in\mathbb{R}^{n\times r}$, where $r\ll n$. Then 
in the \emph{online} phase where we solve many more 
linear systems $A(p_*)x(p_*)=b(p_*)$ for many new values of $p_*$, the idea is to find an approximate solution 
in the span of $X$, that is, $x(p_*)\approx \hat x(p_*)=Xc(p_*)$ for a short vector $c(p_*)\in\mathbb{R}^r$. 
Below we take $p_*$ to be fixed to a particular value in $\Omega$; the overarching goal is to solve $A(p_*)x(p_*)=b(p_*)$ for many different values $p_*$. 

A natural idea is then to find $\hat x(p_*)=Xc(p_*)$ that minimises the residual. 
This leads to the least-squares problem 
\begin{equation}
  \label{eq:snapshotLS}
\min_{c(p_*)}\|A(p_*)Xc(p_*)- b(p_*)\|_2.
\end{equation}
Minimising the residual in a subspace is a classical idea, and employed e.g. in GMRES~\cite{gmres} when the subspace is the Krylov subspace. As we describe in Section~\ref{sec:related}, algorithms based on~\eqref{eq:snapshotLS} (which is 'ApSnap' in our terminology; see~\eqref{eq:diagram}) have been studied in~\cite{markovinovic2006accelerating,guido2024subspace}. 

Note that~\eqref{eq:snapshotLS} is naturally a very tall-skinny least-squares problem, where the coefficient matrix $B:=A(p_*)X$ has size $n\times r$ with $n\gg r$. 
We visualise this using a diagram, where we drop $(p_*)$ for simplicity:
\begin{equation}
  \label{eq:fullsamplepic}  
\min_c\left\|\    \begin{tikzpicture}[node distance=1mm,baseline=-1mm]
       \node[draw, rectangle, solid, fill=gray!10, minimum width=\rrowsize*1.5, minimum height=\colsize] (B) {$\Scale[1]{B}$};
   \end{tikzpicture}\hspace{0.3mm}
   \begin{tikzpicture}
[node distance=1mm,baseline=-1mm]
       \node[draw, rectangle, solid, fill=gray!10, minimum width = 0.3cm, minimum height=\rowsize] (x) {$\Scale[.8]{c}$};
   \end{tikzpicture}
-  \begin{tikzpicture}[node distance=1mm,baseline=-1mm]
       \node[draw, rectangle, solid, fill=gray!10, minimum width=0.3cm, minimum height=\colsize] (b) {$\Scale[.8]{b}$};
   \end{tikzpicture}
\ \right\|_2.
\end{equation}

The key idea of this work is to solve~\eqref{eq:snapshotLS}, or equivalently~\eqref{eq:fullsamplepic},approximately and efficiently by \emph{subsampling} the rows. That is, 
instead of $\min_c\|Bc-b\|$ we solve 
$\min_{\hat c}\|S(B\hat c-b)\|$, where $S\in\mathbb{R}^{s\times n} (r\leq s\ll n)$ is a fat, subsamping matrix, i.e., $S$ is a row-submatrix of $I_n$, each row having exactly one $1$, and $0$ elsewhere. $SB$ then becomes a submatrix of $B$ consisting of a subset of the rows: 
\begin{equation}  \label{eq:subsamplepic}
\min_c\left\|\    \begin{tikzpicture}[node distance=1mm,baseline=-1mm]
       \node[draw, rectangle, solid, fill=gray!10, minimum width=\rrowsize*1.5, minimum height=\colsize] (B) {$\Scale[1]{B}$};
       \node[draw, rectangle, solid, fill=blue!20, minimum width=\rrowsize*1.5, minimum height=\colsize/50, right=-88mm of B, above=-16mm of B] (BB) {};
       \node[draw, rectangle, solid, fill=blue!20, minimum width=\rrowsize*1.5, minimum height=\colsize/50, right=-88mm of B, above=-4mm of B] (BB) {};
       \node[draw, rectangle, solid, fill=blue!20, minimum width=\rrowsize*1.5, minimum height=\colsize/50, right=-88mm of B, above=-8mm of B] (BB) {};
       \node[draw, rectangle, solid, fill=blue!20, minimum width=\rrowsize*1.5, minimum height=\colsize/50, right=-88mm of B, below=-3mm of B] (BB) {};
   \end{tikzpicture}\hspace{0.3mm}
   \begin{tikzpicture}
[node distance=1mm,baseline=-1mm]
       \node[draw, rectangle, solid, fill=gray!10, minimum width = 0.3cm, minimum height=\rowsize] (x) {$\Scale[.8]{c}$};
   \end{tikzpicture}
-  
\begin{tikzpicture}[node distance=1mm,baseline=-1mm]
       \node[draw, rectangle, solid, fill=gray!10, minimum width=0.38cm, minimum height=\colsize] (b) {$\Scale[.8]{b}$};
       \node[draw, rectangle, solid, fill=green!20, minimum width=0.38cm, minimum height=\colsize/50, right=-88mm of b, above=-16mm of b] (bb) {};
       \node[draw, rectangle, solid, fill=green!20, minimum width=0.38cm, minimum height=\colsize/50, right=-88mm of b, above=-4mm of b] (bb) {};
       \node[draw, rectangle, solid, fill=green!20, minimum width=0.38cm, minimum height=\colsize/50, right=-88mm of b, above=-8mm of b] (bb) {};
       \node[draw, rectangle, solid, fill=green!20, minimum width=0.38cm, minimum height=\colsize/50, right=-88mm of b, below=-3mm of b] (bb) {};
   \end{tikzpicture}
\ \right\|_2
\rightarrow \quad 
\min_{\hat c}\left\|\    \begin{tikzpicture}[node distance=1mm,baseline=-1mm]
       \node[draw, rectangle, solid, fill=blue!20, minimum width=\rrowsize, minimum height=\colsize/2] (SB) {$\Scale[1]{SB}$};
   \end{tikzpicture}\hspace{0.3mm}
   \begin{tikzpicture}
[node distance=1mm,baseline=-1mm]
       \node[draw, rectangle, solid, fill=gray!10, minimum width = 0.3cm, minimum height=\rowsize] (x) {$\Scale[.8]{\hat c}$};
   \end{tikzpicture}
-  \begin{tikzpicture}[node distance=1mm,baseline=-1mm]
       \node[draw, rectangle, solid, fill=green!20, minimum width=0.3cm, minimum height=\colsize/2] (Sb) {$\Scale[.6]{Sb}$};
   \end{tikzpicture}
\ \right\|_2.   
\end{equation}

We call the algorithm SubApSnap, 
following the structure of the coefficient matrix $SB = SA(p)X$:
\begin{equation}\label{eq:diagram}
SB = \underbrace{\begin{tikzpicture}[node distance=1mm,baseline=-1mm]
       \node[draw, rectangle, solid, fill=blue!20, minimum width=2*\bigsize, minimum height=\colsize/2.5] (S) {$\Scale[1]{S}$};
   \end{tikzpicture}}_{
   {{\rm {\bf \rr{Sub}}sample}}\atop{s\times n}
}
\hspace{0.3mm}
   \underbrace{\begin{tikzpicture}[node distance=1mm,baseline=-1mm]
       \node[draw, rectangle, solid, fill=gray!20, minimum width=2*\bigsize, minimum height=2*\bigsize] (S) {$\Scale[2]{A(p)}$};
   \end{tikzpicture}}_{
   {{\rm \bf \rr{Ap}}}\atop{n\times n}}
   \hspace{0.3mm}
   \underbrace{\begin{tikzpicture}[node distance=1mm,baseline=-1mm]
       \node[draw, rectangle, solid, fill=red!20, minimum width=\smallsize, minimum height=2*\bigsize] (S) {$\Scale[.8]{X}$};
   \end{tikzpicture}}_{
   {{\rm {\bf \rr{Snap}}shot}}\atop{ n\times r}
   },
\end{equation}
where we reiterate that $n\gg s\geq r$. When $s=r$, we have $SB\hat x = Sb$, and $A(p_*)\hat x$ is equal to $b$ (interpolated) in the $r$ indices selected by $S$. We assume that $SA(p)X$ is full column rank, which holds generically and essentially always when $S$ is chosen sensible, including the methods described in Section~\ref{sec:subsamplechoice}. 

Importantly, we use the same $S$ for all new values of $p_*$, at least those close to a snapshot point $p_i$ from which $S$ is obtained (see Section~\ref{sec:variants}). 
As a result, the SubApSnap solution can be computed very efficiently: 
as $\min_{\hat c}\|S(B\hat c-b)\|_2$ 
is a least-squares problem of much smaller size $s\times r$ where $s=O(r)$, solving it requires only $O(r^3)$ operations. 
We also need to compute $SA(p_*)X$, which would cost $O(snr)=O(nr^2)$ if $SA(p_*)$ is dense (the cost is lower when $A(p)$ is sparse). 
Here, there is potentially an equally significant benefit: 
$SB=SA(p_*)X$, computed as $(SA(p_*))X$, 
only requires us to compute a submatrix of $A(p_*)$, i.e. the row submatrix corresponding to those selected by $S$, i.e., there is no need to even look at the whole matrix $A(p_*)$ to get the SubApSnap solution; in fact we only require a fraction $s/n\ll 1$ of them. This can lead to speedup of a factor of up to $n/s =O(n/r)$, in applications where building $A$ at a new parameter $p_*$ requires substantial computation. 

\subsection{Implementation}
We present SubApSnap in pseudocode in Algorithm~\ref{alg:subapsnap}. 
\begin{algorithm}[t]
  \caption{SubApSnap for solving $A(p)x(p)=b(p)$
  for many parameter values $p\in\Omega$.}\label{alg:subapsnap}
  \begin{algorithmic}[1]
  \vspace{0.5pc}
	\Require{Snapshot points $p_1,\ldots,p_r\in\Omega$, 
    new parameter value(s) $p_*\in\Omega$, subsampling strategy (e.g. leverage score sampling).}
    \Ensure{
    Approximate solution for $A(p_*)x(p_*)=b(p_*)$.
    }
\vspace{0.5pc}
	\State	Solve $A(p_i)x(p_i)=b(p_i)$ for $i=1,\ldots,r$, fully or approximately, and set $X=[x(p_1),x(p_2),\ldots,x(p_r)]$. 
    \State (Optional, recommended when $r\gg 1$) Find (approximate and orthonormal) range $Q$ of $X$ using (randomised) SVD, and replace $X\leftarrow Q$. 
    \State Use subsampling strategy on $A(p_i)X$ %
    for some $p_i$ (e.g. median of $\{p_i\}_{i=1}^r$)
    to obtain 
    subsampling matrix $S\in\mathbb{R}^{s\times n}$, where $s=O(r)$. 
    \State 
    For each value of $p_*$ not among the snapshot points, 
    compute $\tilde A = (SA(p_*))X$ and $\tilde b = Sb(p_*)$, 
    solve 
    $\min_c\|\tilde Ac-\tilde b\|_2$, 
    and output $Xc$ as the solution for $A(p_*)x=b(p_*)$. 
\end{algorithmic}
\end{algorithm}

We will discuss in detail the snapshot and subsampling strategies in the next sections. 
Variants of SubApSnap are discussed in Section~\ref{sec:variants} after we cover some theory. 
Here we discuss other implementation aspects to optimise the speed and numerical stability. 
\paragraph{Computing $SA(p_*)X$.}
As just mentioned, because the SubApSnap least-squares solution~\eqref{eq:subsamplepic} is often so fast, the runtime can be dominated by computing $SA(p_*)X$; in particular, generating the whole $n\times n$ matrix $A(p_*)$ should be avoided. 
It is always advisable to directly compute $SA(p_*)$, that is, the $s$ rows  of $A(p_*)$. This is possible in many applications at much greater speed than forming the whole matrix $A(p_*)$, including those  described in Section~\ref{sec:exp}. 
Another point is that, since $SA(p_i)X$ has been computed at the snapshot points $p_i$, it is sometimes not even necessary to compute $SA(p_*)X$, if $A(p)$ is varying in a structured fashion. For example, in the MOR problem discussed in Section~\ref{sec:mor}, $A(p)=pE-A_0$, so $SA(p_*)X$ can be computed as 
$SA(p_*)X=SA(p_i)X +(p_*-p_i)SEX$, where the first term is known and the second term is simply a submatrix of the snapshot $X$. Note furthermore that once $SEX$ is computed, 
$SA(p)X$ can be computed for another value of $p$ at the mere cost of $O(sr)$, by simply performing the 
scalar multiplication $(p_*-p_i)(SEX)$ and the addition 
$SA(p_i)X +(p_*-p_i)SEX$. More generally, exactly how best to compute $SA(p_i)X$ is problem dependent, and often requires some care; we revisit this in Section~\ref{sec:largedelayMOR}. 

Clearly, the same comment applies to $b(p_*)$, namely, only $s$ of its entries are needed for SubApSnap, and they can often be obtianed via the information on $b(p)$. 
Notably, this cost saving in the evaluation of $b(p_*)$ is what makes DEIM so successful, as we discuss further in Section~\ref{sec:subapsnapvsDEIM}. 

\paragraph{Postprocessing the snapshot matrix $X$.}
We can orthogonalise $X$ by taking its thin QR decomposition $X=Q_XR_X$ and replacing $X$ with $Q_X$; while this is not necessary for the algorithm (the output remains the exact same), it usually results in a better-conditioned problem and thus better accuracy, and it will be convenient to do so in the analysis in Section~\ref{sec:ptheory}. 
In Proper Orthogonal Decomposition~\cite{berkooz1993proper,chatterjee2000introduction} (POD) 
an SVD of $X$ is computed, and $X$ is replaced with its leading left singular vectors (which spans its range). This is particularly effective when $r\gg 1$ and the snapshot matrix $X$ is numerically rank deficient. 
Moreover, recent work~\cite{guido2024subspace} shows that one may benefit from using a randomised rangefinder~\cite{halko2011finding} to efficiently find the range, which we revisit in Section~\ref{sec:related}. All these can be incorporated in SubApSnap.


\subsection{Computational complexity}
We summarise the computational cost in Tables~\ref{tab:comlexity} and~\ref{tab:complexity_newp}. 
For simplicity here we assume $A(p)$ is a dense $n\times n$ matrix, and $s=O(r)$ where $r$ is the number of snapshot points, $r\ll n$ and often $r=O(1)$ suffices, say $r=10$. 

\begin{table}[ht]
  \centering
  \begin{minipage}{0.45\textwidth}
    \centering
\caption{Cost for 
  solving $A(p_i)x(p_i)=b(p_i)$ for 
  $\ell$ parameter values $p_1,\ldots,p_\ell$.}
  \label{tab:comlexity}
  \begin{tabular}{c|c}
  & Complexity   \\\hline
  Full & $O(n^3 \ell )$  \\
  ApSnap & $O(n^3 r+  n^2 \ell r)$ \\
  SubApSnap & $O( n^3 r + n \ell r^2 )$ \\
  \end{tabular}
  \end{minipage}\hfill
  \begin{minipage}{0.45\textwidth}
    \centering
  \caption{Cost for each new value of $p$, 
  }
  \label{tab:complexity_newp}
  \begin{tabular}{c|c}
  & Complexity   \\\hline
  Full & $O(n^3)$  \\
  ApSnap & $O(n^2r)$ \\
  SubApSnap & $O(nr^2)$ \\
  \end{tabular}
\end{minipage}
\end{table}

\ignore{
\begin{table}[htbp]
  \centering
  \caption{Cost for 
  solving $A(p_i)x(p_i)=b(p_i)$ for 
  $\ell$ parameter values $p_1,\ldots,p_\ell$.}
  \label{tab:comlexity}
  \begin{tabular}{c|c}
  & Complexity   \\\hline
  Full & $O(n^3 \ell )$  \\
  ApSnap & $O(n^3 r+  n^2 \ell r)$ \\
  SubApSnap & $O( n^3 r + n \ell r^2)$ \\
  \end{tabular}
\end{table}

\begin{table}[htbp]
  \centering
  \caption{Cost for each new value of $p$, 
  }
  \label{tab:complexity_newp}
  \begin{tabular}{c|c}
  & Complexity   \\\hline
  Full & $O(n^3)$  \\
  ApSnap & $O(n^2r)$ \\
  SubApSnap & $O(nr^2)$ \\
  \end{tabular}
\end{table}
}
The dominant cost of SubApSnap is in computing the snapshot $O(n^3r)$, plus solving~\eqref{eq:problem} at $\ell$ values of $p$ with SubApSnap, each requiring 
the computation of $(SA(p_*))X$, which is $O(nr^2)$, assuming a row of $A(p_*)$ can be computed in $O(n)$ or at most $O(nr)$ operations. 

As hinted above, SubApSnap can be even faster than the tables indicate: for example, it is $O(r^3)$ per new value of $p$ when $A(p)$ is sparse with $O(1)$ nonzeros per row, as then computing $A(p_*)X$ can be done in $O(sr^2)=O(r^3)$ operations. 

In case the quality of the SubApSnap solution $\hat x(p_*)$
is not satisfactory, one can invoke another algorithm, e.g. GMRES, with $\hat x(p_*)$ as the initial guess, to improve the accuracy. This has been suggested in e.g.~\cite{guido2024subspace}. Clearly, this immediately increases the complexity to at least $O(n^2)$, so in this paper we focus on SubApSnap without such refinement. 

\subsection{Related methods}\label{sec:related}
The first reference that described what is essentially 
ApSnap~\eqref{eq:snapshotLS} appears to be
Markovinovic and Jansen~\cite{markovinovic2006accelerating}. A
recent paper by Guido, Kressner, and Ricci~\cite{guido2024subspace} uses randomisation to efficiently capture the snapshot subspace by multiplying by a (Gaussian) random matrix; one can view their method as ApSnapG, wherein instead of $
\min_{c(p_*)}\|A(p_*)Xc(p_*)- b(p_*)\|_2$ as in~\eqref{eq:snapshotLS}, one solves $\min_{c(p_*)}\|A(p_*)XGc(p_*)- b(p_*)\|_2$, where $G$  is a $r\times \hat r$ Gaussian matrix with $\hat r<r$. 
We can certainly adopt this approach, which would result in what might be called SubApSnapG; we do not pursue this, to keep the main focus on the key idea of subsampling. In many of the examples we describe, $r$ is small enough, and/or $X$ is numerically full rank, even though it has decaying singular values; so that there is little room for $G$ to help. 

Combining subspace methods with (randomised) sketching for $Ax=b$ 
was introduced in~\cite{nakatsukasa_tropp}; subsampling is not emphasized there as there is no parameter dependence, so there is little to be gained as compared with other sketches, as finding an appropriate subsampling is as expensive as applying a general efficient sketching. 

A key idea in SubApSnap is that the row-indices chosen by the 
subsampling matrix $S$ should be a good choice for many (possibly all) parameter values. This is an idea put forward in~\cite{park2025low}, where the goal is to find a CUR approximation to parameter-dependent matrices. 


\subsection{SubApSnap vs. DEIM(=SubSnap) and ApSnap}\label{sec:subapsnapvsDEIM}
SubApSnap is related to DEIM, a popular method for nonlinear model order reduction~\cite{chaturantabut2010nonlinear}. In particular, they both employ a combination of subsampling and snapshot. 
However, they are not the same: 
In DEIM, the goal is to approximate a high-dimensional vector $y(p)\in\mathbb{R}^n$ depending nonlinearly on $p\in\Omega$. 
It forms a snapshot based on its evaluation\footnote{We might hence call the DEIM snapshot $Y$, but we stick with $X$ here to crytallise the contrast~\eqref{eq:deim} vs. \eqref{eq:subapsnapvsdeim}. 
In DEIM the snapshot is then usually combined with PCA/POD and orthonormalised,  and denoted by $U$.}
$X=[y(p_1),\ldots,y(p_r)]$.
Then DEIM finds a subsampling matrix $S$ using $X$, and solves the subsampled least-squares problem:
\begin{equation}\label{eq:deim}
\min_c\|Xc-y(p_*)\|_2
\quad \mbox{(Snap)}
\qquad \rightarrow \qquad \min_c\|S(Xc-y(p_*))\|_2\quad \mbox{(SubSnap=DEIM)}.    
\end{equation}

By contrast, in SubApSnap, the problem 
we want to solve after the snapshot $X$ is chosen is ApSnap~\eqref{eq:snapshotLS}, 
which we solve via subsampling
\begin{equation}\label{eq:subapsnapvsdeim}
\min_c\|A(p_*)Xc-b(p_*)\|_2\quad \mbox{(ApSnap)}\quad \rightarrow 
\quad \min_c\|S(A(p_*)Xc-b(p_*))\|_2 \quad \mbox{(SubApSnap)}.     
\end{equation}

As we hope is clear by now, DEIM can thus be viewed as SubSnap in our terminology, that is, a special case of SubApSnap where $A(p)=I$ is not present. 

DEIM benefits mainly from not having to sample $y(p_*)$ in all coordinates, only $O(r)$ of them, namely $Sy(p_*)$. 
In the classical DEIM algorithm, $s=r$. Oversampling $s>r$ has been explored in e.g.~\cite{peherstorfer2020stability}. 
In SubApSnap, in addition to only needing $s$ entries of $b$, we 
 only need $s$ rows of $A(p)$. 

To summarise, the main difference between DEIM and SubApSnap is that in DEIM only the right-hand side changes with $p$, whereas in SubApSnap the matrix changes too. 



\section{Snapshot: theory and practice}\label{sec:snap}
The success of SubApSnap clearly relies on the assumption that 
$\mbox{Span}(X)$ captures the solution for all values of $p$ of interest\footnote{One can of course use different snapshots $X$ for different regions of $p$. We do not pursue this to keep the discussion simple.}. Informally, it is perhaps not surprising that $x(p)$ depends smoothly on $p$ when $A(p),b(p)$ do too; except when $A(p)$ goes through a singular matrix, so if sufficiently many snapshots are taken, $\mbox{span}(X)$ would provide a good subspace for all relevant $p$. 

Theoretically, the approximability of a family of solutions in a subspace is governed by the so-called Kolomogov $n$-width~\cite{pinkus2012n}, and for certain classes of problems it has been proven that the $n$-width decays rapidly~\cite{cohen2016kolmogorov,greif2019decay}, therefore, subspace methods are successful. From a practical side, the snapshot matrix $X$ is often observed to have rapidly decaying singular values in many problems beyond such classes. The speed and nature of decay depend on how $A(p)$ and $b(p)$ vary with $p$, and the location of snapshot points $\{p_i\}_{i=1}^r$. 
This qualitatively implies that there is a good approximate solution $x(p_*)=Xc(p_*)$ such that $A(p_*)x(p_*)\approx b(p_*)$, hence $\|A(p_*)Xc(p_*)- b(p_*)\|_2=\delta$ is small.

When specialised to parameter-dependent linear systems, the singular value decay of a snapshot matrix is proved by Kressner and Tobler~\cite[\S~2]{kressner2011low} under natural smoothness/analyticity assumptions on the dependence of $A(p),b(p)$ on $p$. 
In the more challenging context of extrapolating to future parameter values, (e.g. $p>1$ when $p_i \in [-1,1]$), the existence of a good solution $x(p_*)\in\mbox{Span}(X)$ is shown 
in a recent paper by Guido, Kressner, and Ricci~\cite{guido2024subspace}. More generally, computing a snapshot and working in its dominant subspace is the basis of POD (proper orthogonal decomposition), which is used in a variety of applications, e.g. in control~\cite{chatterjee2000introduction}, fluid dynamics~\cite{berkooz1993proper}, and data science~\cite{brunton2022data}. 

These discussions do not directly guide us on how to choose the snapshot points, when it comes to implementing SubApSnap. This is a subtle issue, and a number of approaches have been proposed, including greedy methods~\cite[\S~3.4]{benner2015survey} and optimisation-based methods~\cite{kunisch2010optimal}, along with grid-based or uniform sampling. 
 As an example to illustrate the complicated nature of the problem, while the analysis in~\cite[\S~2]{kressner2011low} naturally suggests the use of Chebyshev points for snapshot locations, there is usually no noticeable difference between Chebyshev points and equispaced points. On the positive side, we often see that simple approaches are sufficiently effective, for example equispaced or randomly selected snapshot points in $\Omega$. Further, a practical guidance can be gained from the decay of the singular values of the snapshot matrix; sufficient decay usually indicates that enough snapshot points have been taken, assuming $p_i$ has been taken reasonably evenly in the domain $\Omega$. In most of our experiments we took equispaced snapshot points. 

As this is a topic of ongoing research, particularly when $\Omega$ is a high-dimensional set, and because the main novelty in SubApSnap lies in the subsampling, in what follows we \emph{assume} that the problem~\eqref{eq:problem} is amenable to snapshot-based solutions, and that 
the snapshot matrix $X$ gives a good solution at the new parameter value $p_*$ such that there exists $Xc(p_*)$ such that the ApSnap residual
$\|A(p_*)Xc(p_*)- b(p_*)\|_2$ in~\eqref{eq:fullsamplepic} is small.


\section{Subsampling: near-optimality and selection}\label{sec:subsample}
We now turn our attention to the effect of subsamplng. 
Empirically, simple experiments reveal the following:
 The solution for the subsampled 
$\min_{\hat c}\|S(B\hat c-b)\|_2$  gives solutions almost as good as ApSnap, in that
$\|B\hat c-b\|_2\leq C\|B c-b\|_2$ for a modest constant $C$; note that $C$ has to be at least 1 by the definition of $c$. 
Combined with the assumption that the snapshot is good enough so that 
$\|B c-b\|_2\leq \delta$, this means 
\begin{equation}
  \label{eq:defC}
\|A(p_*)(X\hat c)-b(p_*)\|_2\leq C\delta  ; 
\end{equation}
 Typical numbers are e.g. $\delta =10^{-10}$ and $C=2$, in which case we obtain an approximate solution $\hat x(p_*)=X\hat c$ with residual $2\times 10^{-10}$, a good enough solution in most engineering applications.

\subsection{
Near-optimality of subsampling: general bounds}\label{sec:theory}
The remaining question therefore becomes: how to explain the success of SubApSnap~\eqref{eq:subsamplepic} relative to ApSnap~\eqref{eq:fullsamplepic}, that is, 
$C=O(1)$ in \eqref{eq:defC}. 
For this, we first present  general results that relate the solution of a least-squares problem to one solved by subsampling. 
Both bounds~\eqref{eq:boundviaA},~\eqref{eq:boundviaAb} can be found in the literature (e.g.~\cite{chaturantabut2010nonlinear} for~\eqref{eq:boundviaA} and~\cite{woodruff2014sketching,MartinssonTroppacta} for~\eqref{eq:boundviaAb}), but it appears that they are seldom (if ever) presented together in a way that compares and contrasts the merits and drawbacks of each. We also present them in terms of the polar decomposition $B=UH$~\cite{Higham2008FM}, rather than the QR decomposition (even though identical results hold with $U\leftarrow Q$), as it will be necessary for the forthcoming analysis.
Here $U\in\mathbb{R}^{n\times r}$ has orthonormal columns and $H\in\mathbb{R}^{r\times r}$ is symmetric positive semidefinite.

\begin{lemma}\label{lem:SLSbasics}
  Let $B\in\mathbb{R}^{n\times r},b\in\mathbb{R}^{n}$ with $n> r$, and 
  let $c,\hat c$ denote the solutions for the least-squares problems 
$\min_c\|Bc-b\|_2$ and $\min_{\hat c}\|S(B\hat c-b)\|_2$, where $S\in\mathbb{R}^{s\times n}$ with $r\leq s\leq n$. 
Then we have $\|Bc-b\|_2\leq \|B\hat c-b\|_2$. 
Moreover, $\|B\hat c-b\|_2$ is bounded from above as 
\begin{equation}
  \label{eq:boundviaA}
\frac{\|B\hat c-b\|_2}{\|Bc-b\|_2}\leq \frac{\|S\|_2}{\sigma_{\min}(SU)}, 
\end{equation}
where 
$B= UH$ is the polar decomposition.  
We also have 
\begin{equation}
  \label{eq:boundviaAb}
\frac{\|B\hat c-b\|_2}{\|Bc-b\|_2}\leq \kappa_2(S\tilde U)
\end{equation}
where $\tilde U\tilde H = [B\ b]\in\mathbb{R}^{n\times (r+1)}$ is the 
polar decomposition, and
 $\kappa_2(S\tilde U)=\frac{\sigma_{\max}(S\tilde U)}{\sigma_{\min}(S\tilde U)}$. 
\end{lemma}
Note that $S$ is not required to be a subsampling matrix. 
When we do restrict to subsampling $S$, as we approach full sampling $s\rightarrow n$ ($S=I$ in the limit), both bounds tend to an equality with $\frac{\|S\|_2}{\sigma_{\min}(SU)},\kappa_2(S\tilde U)\rightarrow 1$. 

{\sc proof.}
The first inequality $\|Bc-b\|_2\leq \|B\hat c-b\|_2$ follows trivially from the fact that $c$ minimises $\|Bc-b\|_2$. 


Write $B=UH=W\Sigma V^T$ as the polar decomposition and thin SVD respectively\footnote{We use $W$ for the matrix of left singular vectors instead of the more common $U$, as we prefer to reserve $U$ for the unitary polar factor, which plays a central role here.}. 
Write $b = UU^Tb + (I - U U^T)b =: b_1 + b_2$.\\
Note that $c=B^\dagger b = V\Sigma^{-1} W^Tb = (V\Sigma^{-1}V^T)(VW^T)b = H^{-1}U^T b$, so  \[\|Bc-b\|_2=\|UH(H^{-1}U^T b)-b_1 - b_2 \|_2 = \|UU^Tb -b_1 - b_2\|_2 = \|b_2\|_2.\]
The solution of $\min_{\hat c}\|S(B\hat c-b)\|_2$ is 
\begin{align*}
\hat c &= (SB)^\dagger (Sb) = (SB)^\dagger S(b_1 + b_2)= (SUH)^\dagger S(UU^Tb +b_2) \\
&= H^{-1}(SU)^\dagger (SU)U^Tb +H^{-1}(SU)^\dagger Sb_2 = H^{-1}U^Tb + H^{-1}(SU)^\dagger Sb_2 .
\end{align*}
Therefore
\begin{align*}
B\hat c = UH(H^{-1}U^Tb + H^{-1}(SU)^\dagger Sb_2 )
=UU^Tb + U(SU)^\dagger Sb_2
=b_1 + U(SU)^\dagger Sb_2,
\end{align*}
so
\begin{align*}
\|B\hat c -b\|_2 &= \|U(SU)^\dagger S b_2 -b_2\|_2 = \|(U(SU)^\dagger S-I)b_2\|_2 \leq \|U(SU)^\dagger S-I\|_2 \|b_2\|_2 \\
&= \|U(SU)^\dagger S\|_2\|b_2\|_2 \qquad \quad \text{ as $U(SU)^\dagger S$ is a projection~\cite{szyld2006many}} \\
&= \|(SU)^\dagger S\|_2\|b_2\|_2 \leq \|(SU)^\dagger \|_2 \|S\|_2 \|b_2\|_2 \\
&= \frac{ \|S\|_2}{\sigma_{\min}(SU)}\|Bc-b\|_2 .
\end{align*} 

The inequality~\eqref{eq:boundviaAb} is well known, e.g.~\cite{kireeva2024randomized,MartinssonTroppacta}.
\ignore{
We next proceed to prove bound~\eqref{eq:boundviaAb}. 
 Let $[B,b]=\tilde Q\tilde R$ be the QR factorisation. We have 
$S[B,b]=(S\tilde Q)\tilde R$. 
We can write $\|Bv-b\|_2 = \|\tilde Qw\|_2$ for some $w\in\mathbb{R}^{n+1}$, and $\|S(Bv-b)\|_2 = \|(S\tilde Q)w\|_2$. 

 It follows that, using the submultiplicative property $\|XYv\|_2\leq \|X\|_2\|Yv\|_2$
\[
\|S(Bv-b)\|_2 = \left\|S[B,b]\begin{bmatrix}v\\-1  \end{bmatrix}\right\|_2
\leq \|S\tilde Q\|_2
\left\|\tilde R\begin{bmatrix}v\\-1  \end{bmatrix}\right\|_2
\]
and also
\[
\|S(Bv-b)\|_2 \geq \sigma_{\min}(S\tilde Q)\|Bv-b\|_2. 
\]
These hold for any vector $v$. 
Since by the definition of $\hat c$ we have 
$\|S(B\hat c-b)\|_2\leq \|S(Bc-b)\|_2$, 
it follows that
\[
\|B\hat c-b\|_2 
\leq 
\frac{1}{\sigma_{\min}(S\tilde Q)}\|S(B c-b)\|_2\leq 
\frac{\sigma_{\max}(S\tilde Q)}{\sigma_{\min}(S\tilde Q)}\|Bc-b\|_2
=\kappa_2(S\tilde Q)\|Bc-b\|_2. 
\]
}
\hfill$\square$ 
\paragraph{Remarks on Lemma~\ref{lem:SLSbasics}.}
In the literature of randomised linear algebra,~\eqref{eq:boundviaAb} is often stated in the context of the \emph{subspace embedding} property of randomised sketching matrices~\cite{kireeva2024randomized,woodruff2014sketching}. In particular, for any subspace embedding $S$ for $\tilde U$, i.e., $\sigma_i(S\tilde U)\in[1-\epsilon,1+\epsilon]$, we have 
$\frac{\|B\hat c-b\|_2}{\|Bc-b\|_2}\leq \kappa_2(S\tilde U) = \frac{1+\epsilon}{1-\epsilon}$. 
For example, $\kappa_2(S\tilde U)\lesssim\frac{\sqrt{s}+\sqrt{r+1}}{\sqrt{s}-\sqrt{r+1}}$ with high probability 
if $S\in\mathbb{R}^{s\times n}$ is Gaussian~\cite{davidson2001local}, regardless of $A$ and $b$; such embeddings are called \emph{oblivious}. 

In SubApSnap, $S$ is restricted to a subsampling matrix, so such convenient bounds do not necessarily hold. In the next subsection we proceed to examine the SubApSnap context. 

The bound $\kappa_2(S\tilde U)=\frac{\sigma_{\max}(S\tilde U)}{\sigma_{\min}(S\tilde U)}$ is usually better than $\frac{1}{\sigma_{\min}(SU)}$, as clearly $\sigma_{\max}(S\tilde U)=\|S\tilde U\|_2\leq \|S\|_2\|\tilde U\|_2=1$, and often $\sigma_{\max}(S\tilde U)\ll 1$. 
On the other hand $\sigma_{\min}(U)\geq \sigma_{\min}(\tilde U)$, but this difference is usually minor. 


Despite giving a poorer bound, the latter $\frac{1}{\sigma_{\min}(SU)}$ is still of significant interest as it does not involve $b$, and hence holds for \emph{any} right-hand side $b$. This is crucial in 
DEIM, where one would fix the snapshot $X=U$ and solve $\min_{c(p)}\|Uc(p)-y(p)\|_2$ via 
$\min_{\hat c(p)}\|S(U\hat c(p)-y(p))\|_2$ 
for a large number of new parameter values $p$ (recall~\eqref{eq:deim}), so no information is assumed available on the nature of the right-hand side $y(p)$, and hence $\tilde U$ cannot be referred to. 

The problem we are considering is somewhat different, because we are solving the time-dependent system $A(p)x(p)=b(p)$ where $A,b$ depend slowly on $p$. It is intuitively reasonable to expect that a subsampling matrix $S$ that is good for say most values of $p_i$ would be good also for most other (possibly all) values of $p\in\Omega$. 
Below we make this precise. 

\subsection{
Near-optimality of subsampling:  varying $p$}
\label{sec:ptheory}

Our next goal is to show that for unseen values $p_*$ of the parameter, 
a bound analogous to \eqref{eq:boundviaA} holds, assuming that $A(p)$ varies smoothly with $p$. 

In the analysis below it becomes important that $X$ is well conditioned, so we assume $X=Q_X$ (its Q factor in $X=Q_XR_X$). We also recall that 
$A(p)X=B(p)$, 
now making the dependence of $B,b$, and $c$ on $p$ explicit.
First, by~\eqref{eq:boundviaA},
we obtain
\begin{equation}\label{eq:boundattildet}
\frac{\|B(p_*)\hat c(p_*)-b(p_*)\|_2}{\|B(p_*)c(p_*)-b(p_*)\|_2}\leq \frac{1}{\sigma_{\min}(SU(p_*))}    , 
\end{equation}
where $A(p_*)Q_X=U(p_*)H(p_*)$ is the polar decomposition (recall that $X=Q_XR_X$ is the QR factorisation of the snapshot). 
The issue with~\eqref{eq:boundattildet} is that $A(p_*)$ is unknown and we prefer not to refer to it. 
To overcome this, we first derive a bound on 
$\|SU(p_*)-SU(p_0)\|_2$ 
where $p_0$ is a snapshot point close to $p_*$ (i.e., $A(p_0)$ is known; the precise value of $p_0$ depends on $p_*$). 
We then use it to bound 
$|\sigma_{\min}(SU(p_*))-\sigma_{\min}(SU(p))|$, and thus obtain a bound on $\frac{\|B(p_*)\hat c(p_*)-b\|_2}{\|B(p_*)c(p_*)-b\|_2}$. 


\begin{theorem} \label{thm:First_bound}
Fix $p_0\in\Omega$ and let $p_*\in\Omega$ be such that $\|A(p_*)-A(p_0)\|_2\leq \delta$, and 
let $c(p_*)={\argmin}_c\|B(p_*)c-b(p_*)\|_2$ 
and $\hat c(p_*)={\argmin}_{\hat c}\|S(B(p_*)\hat c-b(p_*))\|_2$  
where 
$B(p)=A(p)Q_X$ for an $n\times r$ ($n\geq r$) orthonormal $Q_X$. 
Define 
$C = \frac{3}{\sigma_{\min}(A(p_0)Q_X)}$. Provided that
$\sigma_{\min}(SU(p_0))>C\delta$, we have 
\[\frac{\|B(p_*)\hat c(p_*)-b(b_*)\|_2}{\|B(p_*)c(p_*)-b(b_*)\|_2}\leq \frac{1}{\sigma_{\min}(SU(p_0))-C\delta}. 
\]
 where $S\in\mathbb{R}^{s\times n}$ with $s<n$ is a subsampling matrix. 
\end{theorem}


{\sc proof}.
We first prove that $\|SU(p_0) - SU(p_*)\|_2 $ is bounded by $\frac{3\|S\|_2 }{\sigma_{\min}(A(p_0)Q_X)}\|A(p_0) - A(p_*)\|_2 $, using perturbation theory for the polar decomposition. 
We will use~\cite[Thm.~2]{li1995new}, which shows that the difference between the unitary polar factors 
$\|U(p_0) - U(p)\|_2$
of $A(p_0)Q_X$ 
and $A(p_*)Q_X$
is bounded \footnote{This inequality holds in any unitarily invariant norm. We are using the polar decomposition rather than the more common QR decomposition in our analysis because this nonasymptotic bound \cite[Thm.~2]{li1995new} is available for the former. For QR, we are aware only of first-order bounds.}  by 
$(\frac{2}{\sigma_{\min}(A(p_0)Q_X)+\sigma_{\min}(A(p_*)Q_X)} + \frac{1}{\sigma_{\min}(A(p_0)Q_X)})\|A(p_0)Q_X - A(p_*)Q_X\|_2$. This gives
\begin{align*}
\|SU(p_0) - SU(p_*)\|_2 
&\leq \|S\|_2 \|U(p_0) - U(p_*)\|_2 \\
&\leq (\frac{2}{\sigma_{\min}(A(p_0)Q_X)+\sigma_{\min}(A(p_*)Q_X)} + \frac{1}{\sigma_{\min}(A(p_0)Q_X)})\|A(p_0)Q_X - A(p_*)Q_X\|_2  \\
&\leq  (\frac{3}{\sigma_{\min}(A(p_0)Q_X)})\|A(p_0)Q_X - A(p_*)Q_X\|_2 \\
&\leq  (\frac{3}{\sigma_{\min}(A(p_0)Q_X)})\|A(p_0) - A(p_*)\|_2 \|Q_X\|_2 =\frac{3}{\sigma_{\min}(A(p_0)Q_X)}\|A(p_0) - A(p_*)\|_2, 
\end{align*}
where we used the fact $\|S\|_2=\|Q_X\|_2=1$, as they both have orthonormal columns. 

By Weyl's inequality for singular values we  have $\sigma_{\min}(SU(p_*))\geq \sigma_{\min}(SU(p))-\|SU(p_*)-SU(p)\|_2$, so 
\begin{align*}
\frac{1}{\sigma_{\min}(SU(p_*))}
&\leq 
\frac{1}{\sigma_{\min}(SU(p))-\|SU(p_*)-SU(p)\|_2}\\
&\leq
\frac{1}{\sigma_{\min}(SU(p))-\frac{3\|S\|_2 }{\sigma_{\min}(A(p_0)Q_X)}\|A(p_0) - A(p_*)\|_2}    
\end{align*}
where the denominator is positive by the assumption 
$\sigma_{\min}(SU(p))>C\|A(p_0) - A(p_*)\|_2$. 
Substituting this into \eqref{eq:boundattildet} completes the proof. 
\hfill$\square$


We now derive bounds on the residual near-optimality with subsampling that are applicable for all values of $p$. 
\begin{corollary}\label{cor:maincor}
Following the notation of Theorem~\ref{thm:First_bound}, 
suppose that $A(p)$ is Lipschitz continuous in $p$, that is, $\|A(p)-A(\tilde p)\|_2\leq L\|p-\tilde p\|_2$ for all $p,\tilde p\in\Omega$. 
Then for every $p_*\in\Omega$ we have
\begin{equation}\label{eqn:boundclosest}
    \frac{\|B(p_*)\hat{c}(p_*)-b(p_*)\|_2}{\|B(p_*)c(p_*)-b(p_*)\|_2}\leq \frac{1}{\sigma_{\min}(SU(p_{i\star}))-\widetilde C\widetilde\delta}, 
\end{equation}
provided that $\sigma_{\min}(SU(p_{i\star})) \geq \widetilde C\widetilde\delta$, where 
$p_{i\star}$ is the snapshot point nearest to $p_*$, 
$\widetilde\delta = L\|p_*-p_{i\star}\|$, and $\widetilde C= \frac{3}{\sigma_{\min}(A(p_{i\star})Q_X)}$.  

Suppose further that the snapshot points $p_i$ are taken such that 
$\min_i\|p-p_i\|_2\leq h$ for all $p\in\Omega$\footnote{For example, $h$ will be the grid size when $\Omega\in\mathbb{R}$ is an interval and equispaced snapshots are taken.}. 
Then for every 
$p_*\in\Omega$ 
\begin{equation}\label{eq:bound1}
    \frac{\|B(p_*)\hat{c}(p_*)-b(p_*)\|_2}{\|B(p_*)c(p_*)-b(p_*)\|_2}\leq \frac{1}{\min_i\left(\sigma_{\min}(SU(p_i))\right)-\hat C\hat\delta}
\end{equation}
provided that 
$\min_i\left(\sigma_{\min}(SU(p_{i\star}))\right)>\hat C\hat\delta$. 
Here, $\hat C= \frac{3}{\min_i\left(\sigma_{\min}(A(p_i)Q_X)\right)}$ and $\hat \delta = Lh$. 
\end{corollary}

{\sc proof}.
To prove~\eqref{eqn:boundclosest}, 
we consider applying Theorem \ref{thm:First_bound} with  $p_0\leftarrow p_{i\star}$. 
As $A(p)$ is Lipschitz continuous in $p$, we can  choose $\delta = L\|p_*-p_{i\star}\|_2 \geq \|A(p_*)-A(p_{i\star})\|_2 $. 
Applying Theorem \ref{thm:First_bound} then gives~\eqref{eqn:boundclosest}.

For~\eqref{eq:bound1}, 
we follow the same arguments, and now note that 
$\|p_*-p_{i\star}\|_2\leq h$ by assumption, so 
$\|A(p_*)-A(p_{i\star})\|_2\leq L\|p_*- p_{i\star}\|_2\leq hL=\hat\delta$, and
\[\frac{\|B(p_*)\hat c(p_*)-b(p_*)\|_2}{\|B(p_*)c(p_*)-b(p_*)\|_2}\leq \frac{1}{\sigma_{\min}(SU(p_{i\star}))-C\hat\delta},
\] 
We can further bound this to make a global bound. Clearly $C \leq \hat C$ and $\sigma_{\min}(SU(p_{i\star})) \geq \min_i\left(\sigma_{\min}(SU(p_i))\right)$. Since by assumption $\min_i\left(\sigma_{\min}(SU(p_i))\right)>\hat C\delta$, this gives~\eqref{eq:bound1}.
\hfill$\square$

Comparing the two bounds,~\eqref{eqn:boundclosest} tends to be slightly tighter than~\eqref{eq:bound1}, while~\eqref{eq:bound1}
has the advantage that the bound is a constant and does not depend on $p$.
We should note that they are both likely to be severe overestimates in practice. First, as noted above, the bound~\eqref{eq:boundviaAb} indicates why we might expect $1/\sigma_{\min}(SU)$ (which is known to be exponentially large in the worst case with LUPP subsampling~\cite{chaturantabut2010nonlinear}) to be a significant overestimate. Second, the tools used in Theorem~\ref{thm:First_bound}, namely~\cite[Thm.~2]{li1995new} and Weyl's inequality, give the deterministic worst-case bound for the polar factor and singular values. In practice, we  expect these bounds to be overestimates. 
For example, a random perturbation of a rectangular matrix is more likely to \emph{increase} the smallest singular value than to decrease them~\cite{burgisser2010smoothed,stewart1998perturbation}. 
Nonetheless, these are both essentially tight bounds in that examples exist where one cannot improve them (for~\cite[Thm.~2]{li1995new}, by more than a factor 3), so an improved bound may require different assumptions. 



\begin{figure}[htbp]
    \centering
\includegraphics[height=.4\textwidth]{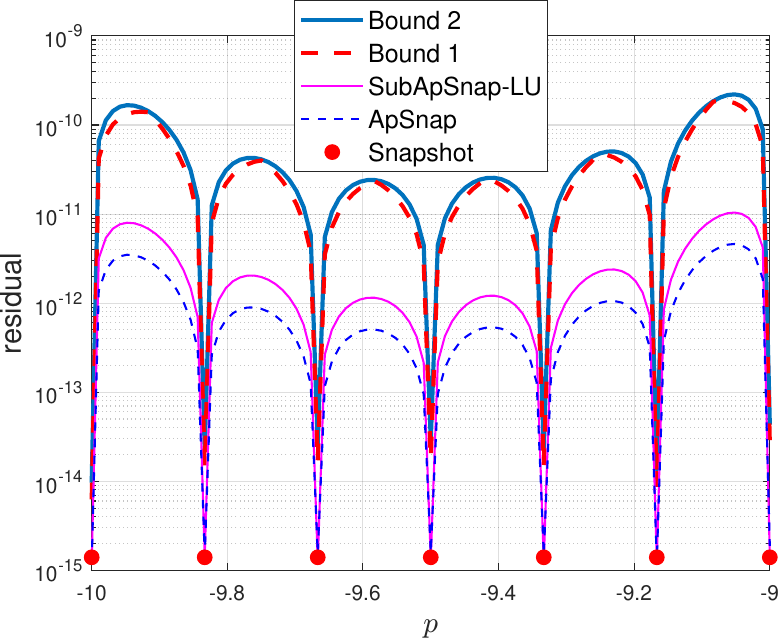}
        \caption{SubApSnap residual and its bounds~\eqref{eqn:boundclosest} (shown as bound 1) and ~\eqref{eq:bound1} (bound 2), for the example~\eqref{ex:tridiag}. The plots are 0 at the snapshot points $p_i$ by construction. 
        }
    \label{fig:boundillustrate}
\end{figure}

To illustrate this, we consider an example where 
\begin{equation}\label{ex:tridiag}
    A(p)=A_0-pI,\qquad b(p) = \exp(b(0)\sin(\frac{p}{10})p), 
\end{equation}
where $A(p)$ is a tridiagonal matrix with 
$A_0$ having $2$ on the diagonals and $-1$ on the first upper and lower diagonal entries, 
 the $\exp$ and $\sin$ in $b(p)$ are taken componentwise, and $b(0)$ is a random Gaussian vector. 
Figure~\ref{fig:boundillustrate} shows the SubApSnap residual and its bound~\eqref{eqn:boundclosest} as $p$ is varied between $[-10,-9]$, together with the ApSnap residuals. 
Note that all the residuals dip to 0 at seven values of $p$; these are the snapshot points $p_i$, where this phenomenon is expected by construction as $X$ contains the solution exactly at the snapshot points. 

Here, the bound, while pessimistic, provides useful information on the quality of SubApSnap solutions; of course, the quality of ApSnap solution $\|B(p_*)c(p_*)-b(p_*)\|_2$ is unknown in practice---in using the above bounds we need to assume this is small enough, and the bounds tell us how close the SubApSnap solution is to ApSnap. 
However, as discussed above, the bound is almost always an overestimate (except for values of $p$ near a snapshot point) and can be inapplicable when the assumption is violated, i.e. the denominator is negative. One can easily construct examples where the bound is not useful. 
Equally, for simple problems the ApSnap/SubApSnap solutions can be so good that the (exact) residuals can easily go below machine precision, and in such cases the bounds are also difficult to use as they do not account for rounding errors. 

The Lipschitz continuity assumption is weaker than e.g. analyticity. It would be interesting to see if stronger bounds can be obtained under stronger assumptions. Equally, the Lipschitz consant $L$ is not always known in applications. Other ways to quantify the SubApSnap accuracy would be desirable. From a practical viewpoint, we describe an estimate of the residual shortly in Section~\ref{sec:levresest}. 

\subsection{Choice of subsampling matrix}\label{sec:subsamplechoice}
Many algorithms are available (e.g.~\cite{dong2023simpler}) for choosing the row subset selection/subsampling matrix $S$ (known as the column subset selection problem): LU with partial pivoting (LUPP) applied to $A(p_i)X$ and column-pivoted QR (CPQR) applied to $(A(p_i)X)^T$ 
are among the classic and well-known approaches. 
Osinsky's recent work~\cite{osinsky2025close} is a recent deterministic algorithm with powerful theoretical guarantees. 
Randomised methods 
that give high-quality guarantee in expectation are also available. 
These include 
leverage score sampling~\cite{drineas2012fast} (sampling proportional to the squared row-norms of $Q$, where $A(p_i)X=QR$), 
and the determinantal point process~\cite{derezinski2021determinantal}, which comes with strong theory but is typically expensive to implement, and 
a recent work by Cortinovis-Kressner~\cite{cortinovis2024adaptive} proposes ARP (adaptive randomised pivoting). 
In LU, QR, Osinsky, and ARP, $s=r$, i.e., interpolation is enforced. 
Alternatively, one can find the pivots via the concatenated matrix $[A(p_i)X\ b]$, which would return $s=r+1$ pivots. In practice, we observe that all these approaches almost always work very well. 
We illustrate this in Figure~\ref{fig:MORsimple} (left panel), where these methods consistently outperform a randomly chosen $S$. The plot shows a trend that is empirically common: randomly chosen pivots are unreliable, ARP is often outperformed by LU and QR---with QR often being slightly more accurate---and leverage score sampling often gives the best accuracy especially when $s$ is sufficiently large ($s=4r$ in all our experiments; the theory $s=O(r\log r)$, and we take $r=O(1)$). 

Regarding choosing the method, when speed or reducing $s$ is of primary importance, it is sensible to use the simplest/fastest method, which is often LUPP~\cite{dong2023simpler}. 
If high accuracy (near-best residual close to ApSnap) is desirable, and a strong theoretical guarantee (in expectation), we recommend using leverage scores, at the expense of oversampling, that is, $s>r$, slightly increasing the sampling cost and the size of the least-squares problem\footnote{When using a method that oversamples $s>n$, including leverage scores, 
we can/should be using 
\emph{weighted} samples, i.e., with leverage-score sampling the $1$'s in $S$ will be replaced with $1/\sqrt{p_i}$ where $p_i$ are the leverage scores, proportional to the probability at which the $i$th row is sampled. Accordingly, left-multiplying $S$ results in subsampling and weighting, and hence $\|S\|_2\neq 1$. The bounds in Section~\ref{sec:theory} accordingly need to be modified. For example, 
the bound in Corollary~\ref{cor:maincor}
should be replaced with $\frac{\|B(p)\hat{c}(p)-b(p)\|_2}{\|B(p)c(p)-b(p)\|_2}\leq \frac{\|S\|_2}{\sigma_{\min}(SU(p_{\text{mid}}))-\hat C\|S\|_2\delta}$. 
}. 

\subsubsection{Residual estimation with leverage score sampling} \label{sec:levresest}
When a modest oversampling is not an issue (as we usually expect to be the case), we recommend a subsampling strategy that oversamples, such as 
leverage score sampling, for the following two reasons: (i) It naturally allows for (or requires) oversampling $s> r$, which tends to result in better accuracy 
with a modest amount of oversampling, i.e., $s=O(r\log r)$. 
(ii) Oversampling, together with the associated weighting, allows us to estimate the (unknown) residual $\|A(p_*)\hat x-b(p_*)\|_2$ via the (known) subsampled version
$\|S(A(p_*)\hat x-b(p_*))\|_2$. 

\begin{figure}
  \begin{minipage}[t]{0.5\hsize}
        \includegraphics[width=0.95\linewidth]{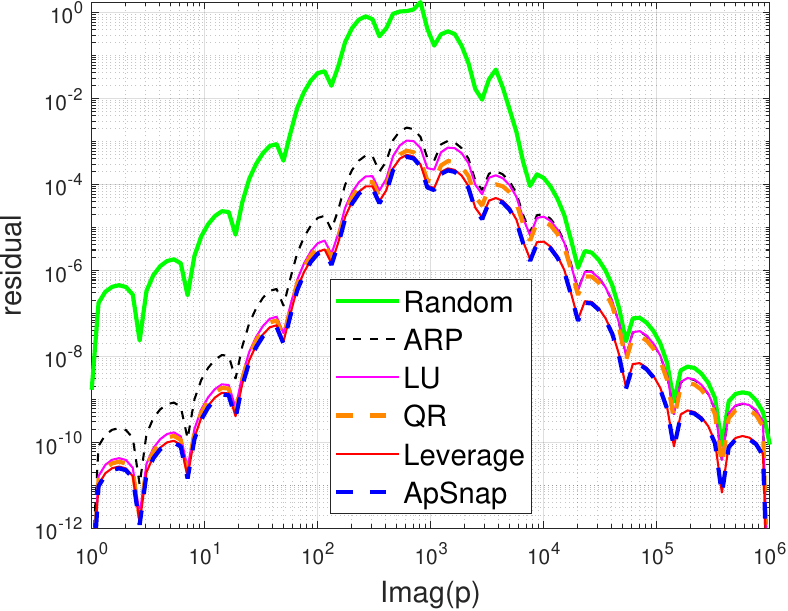}  
  \end{minipage}   
  \begin{minipage}[t]{0.5\hsize}  
          \includegraphics[width=0.95\linewidth]{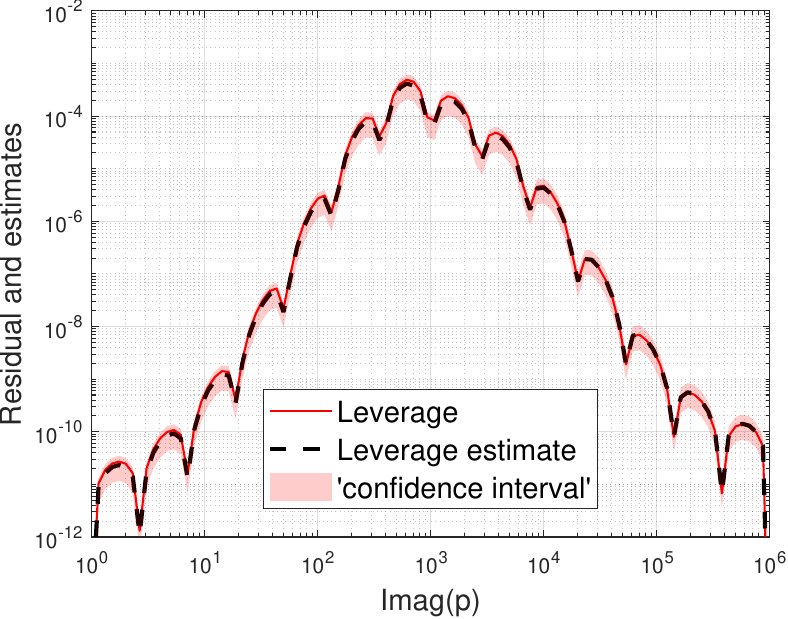}
  \end{minipage}  
    \caption{
    Left: Comparison of subsampling methods for a 
   model order reduction problem; details are given in Section~\ref{sec:mor}. Right: same plot only showing leverage score sampling, together with its estimate and heuristic confidence interval. 
    }
    \label{fig:MORsimple}
\end{figure}

That is, with leverage score sampling on $\mbox{span}([B, b])$, we have with high probability~\cite[\S~9.6.3]{MartinssonTroppacta} 
$(1-\epsilon)\|Bc-b\|_2\leq \|S(Bc-b)\|_2\leq (1+\epsilon)\|Bc-b\|_2$
for $\epsilon \approx  \frac{r\log r}{s}$ (e.g. $\epsilon\approx 1/2$), where $B$ is $n\times r$ as always. 
Note that this gives estimates for the actual (unseen) residual without looking at the whole $B$: $\|Bc-b\|_2\approx \|S(Bc-b)\|_2$, 
together with upper and lower bounds $(1\pm \epsilon) \|S(Bc-b)\|_2$. 
This gives reliable estimates at the $p$-value where the leverage scores are computed. In practice in SubApSnap we use the same subsampling for other values of $p$, so these estimates are only heuristic estimates. Nonetheless, empirically we observe them to give good estimates, providing useful information. 
This is illustrated in the right panel of Figure~\ref{fig:MORsimple}, where
the shaded region indicates the lower and upper estimate of the residual of the SubApSnap solution, which indicates that the estimates capture the actual residual quite accurately. 
Note that such estimates are not possible with interpolation-based methods where $s=r$, as then $S(Bc-b)=0$.

In addition to the ability of estimating the residual, as we have seen in experiments presented and others, leverage score sampling tends to also provide some of the best-quality solutions among the pivoting strategy tested.
Our recommendation therefore is to use leverage score sampling. When oversampling is not desirable, we suggest either LU or QR, as these deterministic approaches tend to work very well in practice. When the randomness in leverage-score sampling is undesirable, the BSS sampling~\cite{batson2009twice} may be a sensible choice.

\subsubsection{Variants of subsampling strategies}\label{sec:variants}

The theory presented in this section suggests a few variants for choosing the subsampling matrix $S$ in SubApSnap. 
First, one can change $S$ depending on the value of $p_*$, that is, choose the $S$ appropriate to its closest snapshot point $p_{i*}$, by pivoting on $A(p_{i*})X$. The rationale is that  this way $S$ is more likely to be a good choice for $p_*$, than one based on a snapshot $p_j$ far away. This comes with $r$ additional runs of the pivoting strategies, which roughly amounts to $O(nr^3)$ operations instead of $O(nr^2)$. 

Another possibility is to chose $r$ indices at several (say $t$) snapshot points, and collect the $rt$ locations (where we expect significant  overlap), and take the union. This would naturally let us oversample at each SubApSnap run. 

These are all effective strategies that can improve the robustness of SubApSnap. 
In our experiments, the simplest strategy of reusing the $S$ chosen at a single representative $p_i$ was sufficiently effective, so we did not pursue these more expensive variants. 


\section{Examples and numerical experiments}\label{sec:exp}
Where speed performance is reported, the experiments were conducted in 
MATLAB version 2025a on a workstation with 512GB RAM. 
Unless otherwise stated, we computed the subsampling matrix $S$ by applying a pivoting strategy (LUPP, leverage scores etc) on $A(p_{i*})X$, where $p_{i*}$ is the median of the snapshot points.



\subsection{Parameter-dependent PDE}
\begin{figure}[htbp]
    \centering
\includegraphics[width=1\linewidth]{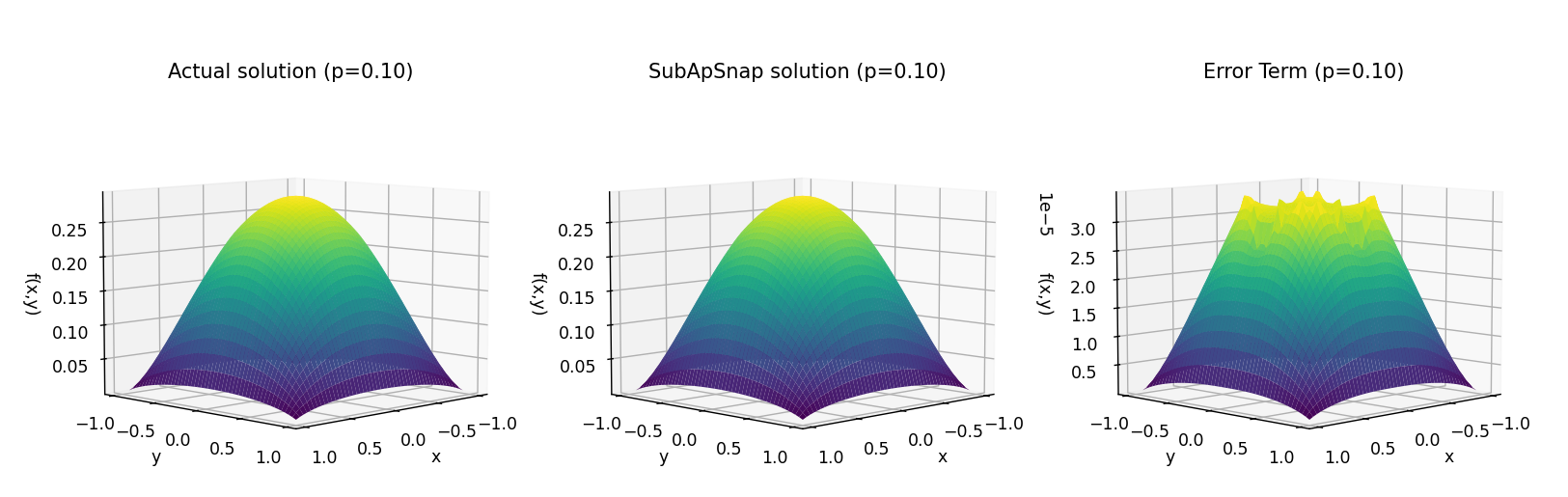}

\includegraphics[width=1\linewidth]{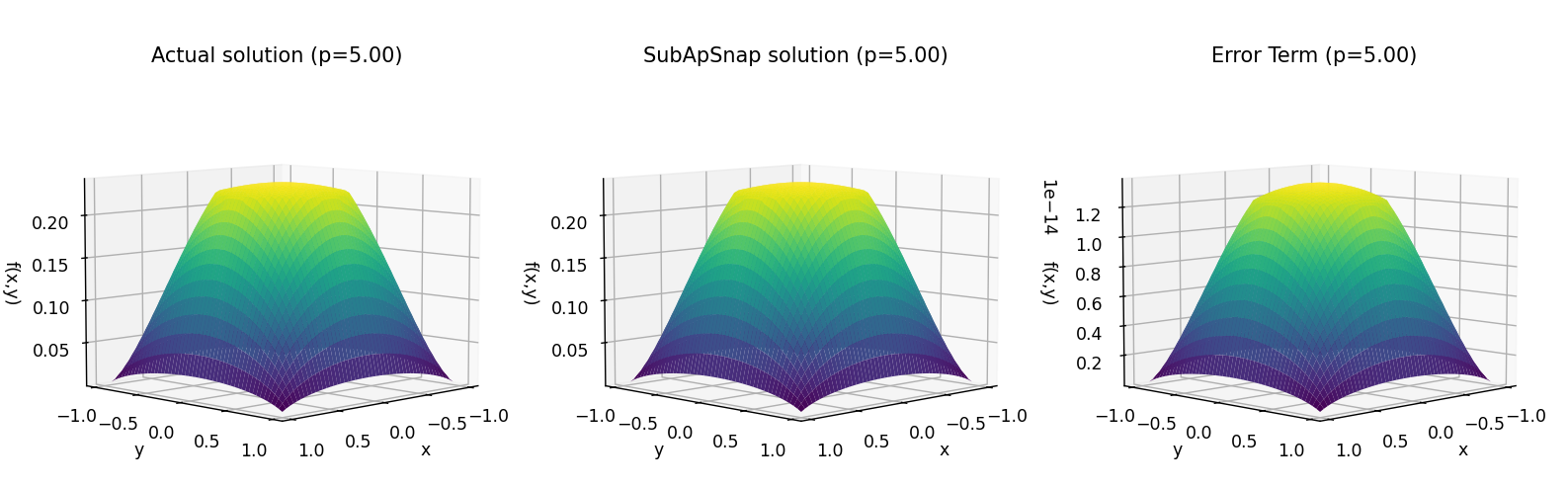}

\caption{Solution and error for 2-d heat equation with SubApSnap for two values of $p$.}
    \label{fig:heat}
\end{figure}

We consider the test problem in~\cite[\S~2.3.1]{kressner2011low}: the stationary heat equation in a 2-d square
$\nabla (\sigma_p(x)\nabla u)=f$ in $[-1,1]^2$, and $u=0$ on the domain boundary. 
Here $\sigma_p(x)$ represents the heat conductivity, and depends on a parameter $p$ as follows:
$\sigma_p(x)=1+p$ in the unit disk $\|x\|_2\leq 1$, and $\sigma_p(x)=1$ otherwise. 
We used the standard finite difference method on a 
on a $100\times 100$ uniform grid 
to discretise the problem
for varying values of $p$, which results in $A(p)x(p)=b$, 
for which we use a standard solver (MATLAB's backslash) and SubApSnap. The solutions 
and error are depicted in Figure~\ref{fig:heat}. With just $r=5$ snapshot points, SubApSnap is able to get solutions of $O(10^{-5})$ accuracy across $p\in [0,5]$. 

While previous approaches have successfully been applied to such model problems, 
SubApSnap is able to deal with the case where the dependence on $p$ is non-smooth, e.g. involving $|p|$. 



\subsection{Transfer function evaluation in model order reduction}\label{sec:mor}
We now describe experiments with SubApSnap for evaluating a 
\emph{transfer function} of the form
\begin{equation}
  \label{eq:transfer_func}
  H(p) = c^T(pE-A)^{-1}b.
\end{equation}
Here, $A,E\in\mathbb{R}^{n\times n}$, and $b,c\in\mathbb{R}^{n}$. 
This results from a linear dynamical system 
\begin{align*}
  E\dot x(t) = Ax(t) + bu(t),
 \qquad y(t) = c^Tx(t),
\end{align*}
and taking the Laplace transform~\cite[Ch.~2]{antoulas2020interpolatory}\footnote{Usually the transfer function takes $s$ as the input variable, but here we use $p$ for consistency with the rest of the paper.}.

We consider evaluating $H(p_*)$ at many values of $p_*$. For each $p_*$, one needs to solve the linear system $(p_*E-A)x(p_*)=b$; thus SubApSnap is applicable. 
 This can be useful in performing model order reduction (MOR), where one attempts to find an approximation $H(p)\approx \tilde H(p) = \tilde c^T(p\tilde E-\tilde  A)^{-1}\tilde  b$
where $\tilde  A$ is $\hat n\times \hat n$ for ideally $\hat n\ll n$; evaluating $H(p)$ at many $p$ is usually required in the MOR algorithm~\cite{antoulas2020interpolatory}. 

We now explain the problem solved in Figure~\ref{fig:MORsimple}. 
We took the convection-diffusion problem in lyapack~\cite{lyapack}, taking $n=2500$. 
The domain of interest is $\Omega = {\rm i}[1,10^6]$. We took $r=15$ snapshot points 
at logarithmically spaced points on $\Omega$. 


\subsubsection{Medium-scale computation}
We now take a larger-scale transfer function 
where 
$A,E\in\mathbb{R}^{n\times n}$, 
$n=66,917$, taken from the Gas Sensor problem from the MOR wiki~\cite{morwiki_gas}.
We sample the transfer function at $5000$ logarithmically spaced points on ${\rm i}[10^{-2},10^4]$, and take $r=30$ snapshots for SubApSnap. 

\begin{figure}[htbp]
    \centering
\includegraphics[height=.4\textwidth]{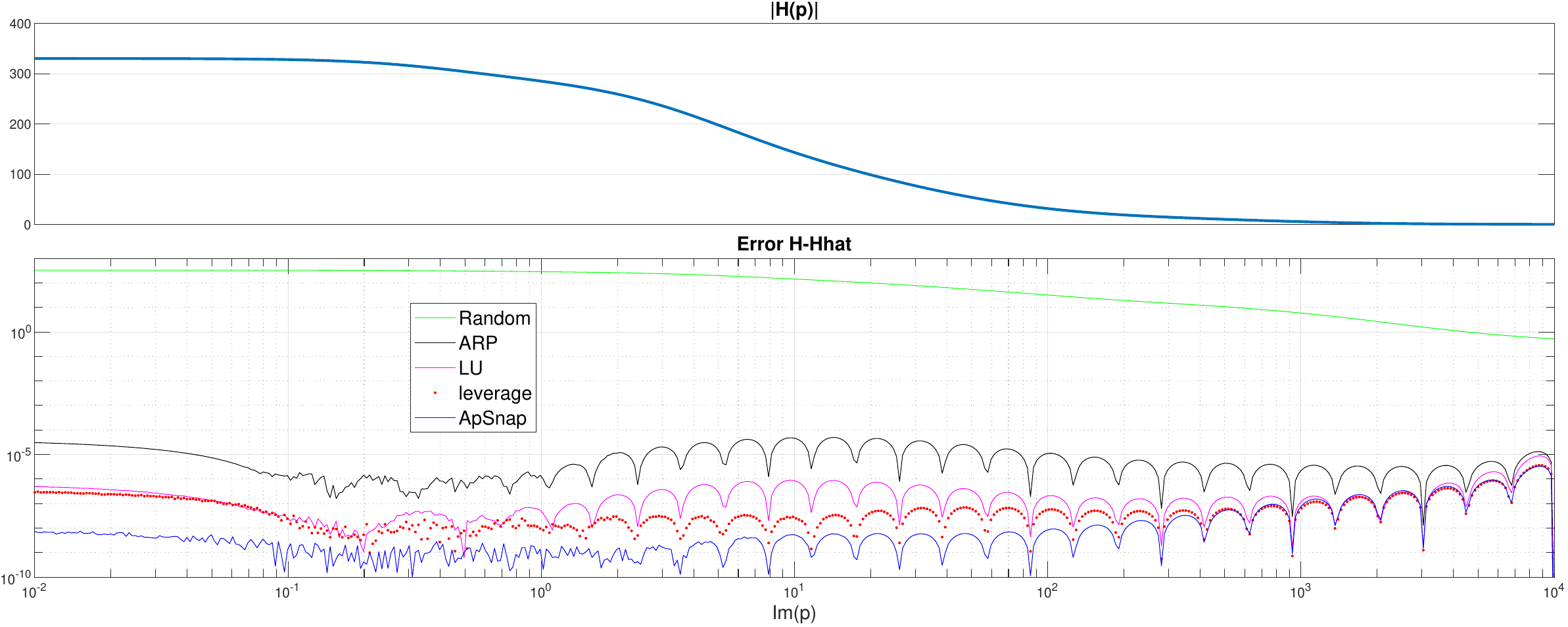}        
    \caption{Transfer function $H(p)$ (above) and errors of their approximations with 
    SubApSnap with four subsampling strategies 
    (Random, LU with partial pivoting, ARP, and leverage-score sampling), and ApSnap (below).
      }
    \label{fig:bigMOR_error}
\end{figure}

Figure~\ref{fig:bigMOR_error} shows the SubApSnap errors 
with different subsampling strategies; note that we are plotting the error in the transfer function $|H(p)-\hat H(p)|$ where $\hat H(p)=c^T\hat x(p)$, not the residual $\|(pE-A)\hat x-b\|_2$. With $r=30$ snapshot points, SubApSnap is able to give five digits of accuracy across $p\in\Omega$; we can further improve accuracy by increasing $r$. 
Among the subsampling strategies, leverage score tends to give the best accuracy, followed by LU (and QR, not shown here). 
Random sampling gives no accuracy with $\tilde H(p)=0$, even at the snapshot points; this is because $b$ is sparse with $744$ nonzeros, and none of the 30 subsamples picked any of them.  

Table~\ref{tab:runtimes} reports the execution time for evaluating $H(p)$ at 5,000 new $p$-values. 
We observe more than $2,500\times $ speedup in the online phase. 
If we include the snapshot time as part of SubApSnap the speedup is about $\times 17$;  this factor of course increases with the number of online evaluations. The table also shows the error in the approximation of the transfer function and the residuals $\|A(p)x(p)-b\|_2/\|b\|_2$, for which the maximum and median values are printed over the $p$-values. 

\begin{table}[htbp]
  \centering
  \caption{Runtime and error for a medium-scale ($n=66,917$) transfer function approximation.}
  \label{tab:runtimes}
  \begin{tabular}{c|c|c|c}
  & time(s) & error $\max_p|H(p)-\hat H(p)|$ & (max, median) $\|A(p)x(p)-b\|_2/\|b\|_2$\\\hline
  Full & 5231.3 & &\\
  Snap & 299.3 &  &\\
  ApSnap & 50.6 & 3.3e-6 & 8.8e-6, \quad 1.2e-9\\
  SubApSnap-LU & 1.9 & 8.7e-6 & 1.3e-5, \quad 3.6e-9\\
 SubApSnap-ARP & 2.0 & 4.9e-5 & 8.5e-5, \quad 5.1e-9\\
  SubApSnap-Lev & 2.2 &  2.3e-6   & 1.0e-5, \quad  1.4e-9
  \end{tabular}
\end{table}

\begin{figure}
  \begin{minipage}[t]{0.5\hsize}
          \includegraphics[width=0.95\linewidth]{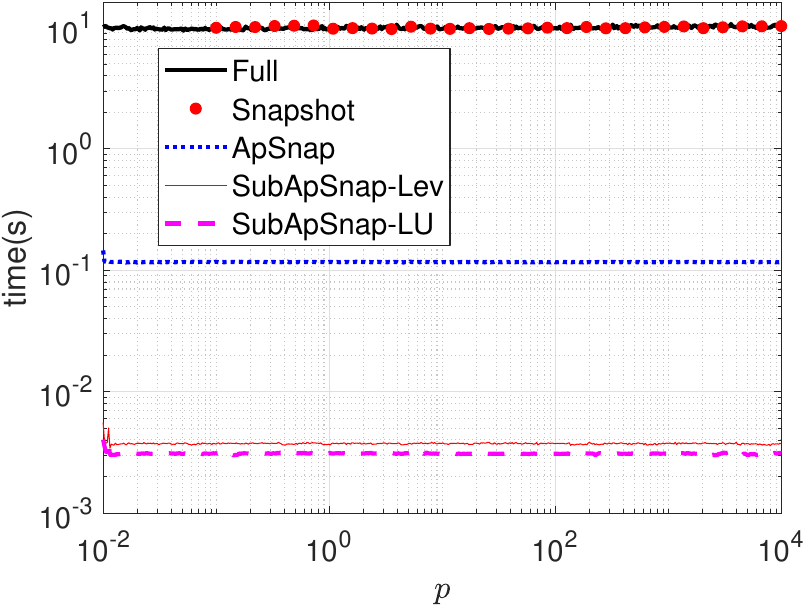}
  \end{minipage}   
  \begin{minipage}[t]{0.5\hsize}
      \includegraphics[width=0.95\linewidth]{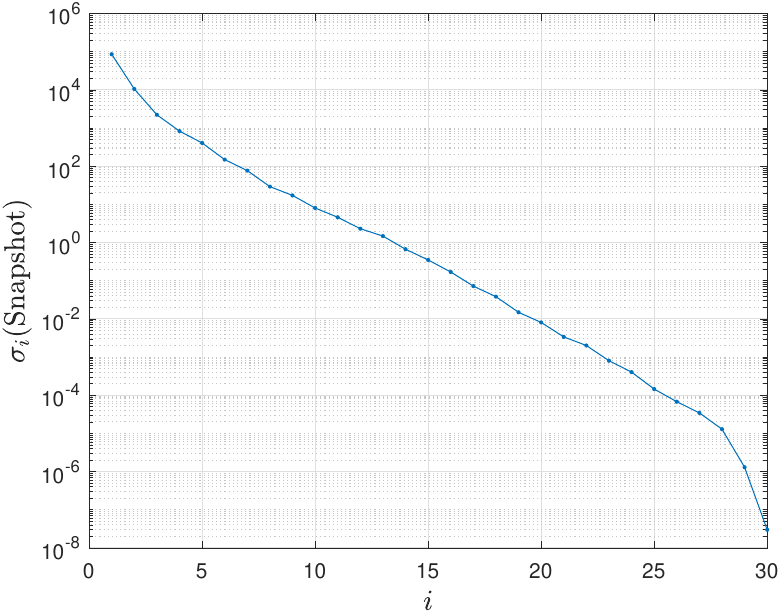}  
  \end{minipage}  
    \caption{
   Model order reduction problem. Left: runtime for different values of  $p$. 
   Note that snapshots are taken only at $r=30$ points $p_i$, whereas other methods are used for $500$ values of $p$. 
   Right: singular values of the snapshot matrix $X$ (before QR or SVD preprocessing). 
    }
    \label{fig:MORsvd}
\end{figure}

Figure~\ref{fig:MORsvd} shows the time required for each linear system as $p_*$ is varied, and the decay of the singular values of the snapshot matrix $X$. 

\subsubsection{Larger and more complicated functions; delay system}\label{sec:largedelayMOR}
 We next follow~\cite{beattie2009interpolatory} and consider a linear dynamical system with an internal delay, which results in the transfer function of the form  $H(p) = c^T(pI-A_0-e^{\tau p}A_1)^{-1}b$, where $\tau>0$, $\kappa>2$ are  constants (we take $\tau=0.1,\kappa=2.1$), and 
 $A_1 = \frac{1}{\tau}(T-\kappa I), A_0=3A_1$, where $T_{ij}=1$ when $|i-j|=1$ and $i=j=1,i=j=n$, otherwise $T_{ij}=0$. Clearly, the main computational task is to solve $A(p)x(p)=b$ where $A(p)=pI-A_0-e^{\tau p}A_1$. 

 We take $n=10^7$, a large-scale example. 
 For each new $p_*$, it is important to compute $SA(p_*)X$ as efficiently as possible in SubApSnap. Here we do so by precomputing 
 $SX, SA_0X$, and $SA_1X$, and take $SA(p_*)X=pSX-SA_0X-e^{\tau p}SA_1X$. In addition, we precompute $c^TX$, and find $H(p_*)\approx (c^TX)\hat c$, where $\hat c$ is the solution for $\min_{\hat c}\|SA(p_*)X\hat c-Sb\|_2$. This way, each offline step  completely avoids any computation involving $n$-dimensional vectors or matrices. 
 
\begin{figure}
    \centering
\includegraphics[width=0.9\linewidth]{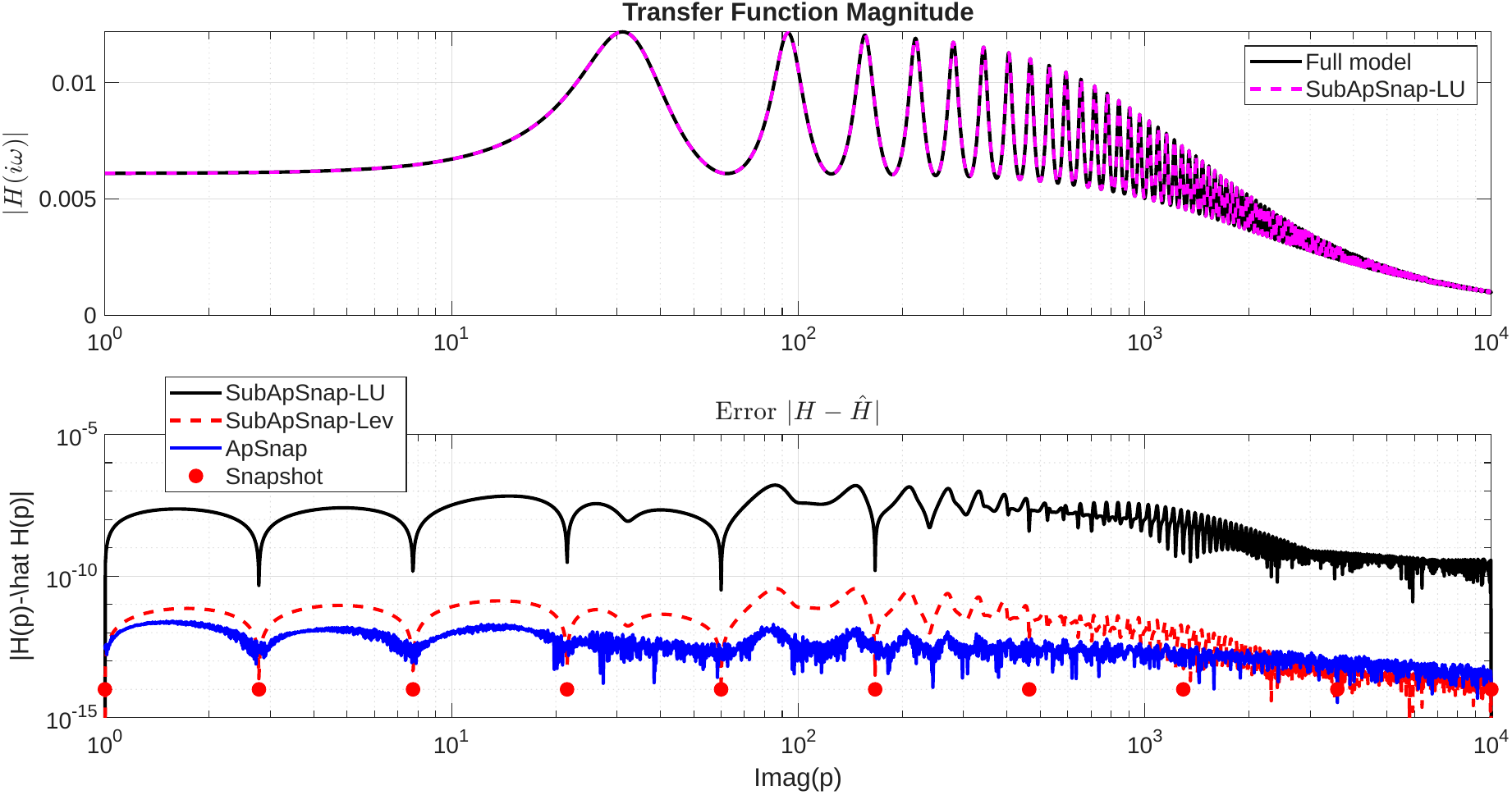}    

    \caption{
    Transfer function and error for the delay system with $n=10^7$. 
    }
    \label{fig:delay}
\end{figure}

 The transfer function along with the SubApSnap error is shown in Figure~\ref{fig:delay}, and Table~\ref{tab:delay} reports the runtime. 
 
 This example shows that even when $A(p)$ is sparse (tridiagonal here, making ApSnap relatively inefficient) and amenable to linear solves, SubApSnap can offer dramatic speedups; it is achieving $30,000\times$ with LU, and $20,000$ with leverage score sampling (which gives more accurate results).


\begin{table}[htbp]
  \centering
  \caption{Runtime and error for a sparse, large-scale ($n=10^7$) transfer function approximation.}
  \label{tab:delay}
  \begin{tabular}{c|c|c|c}
  & time(s) & error $\max_s|H(s)-\hat H(s)|$ & (max, median) $\|A(p)x(p)-b\|_2/\|b\|_2$\\\hline
  Full &    9806.9 & &\\
  Snap &  47.3 &  &\\
  ApSnap & 14386.1 & 2.6e-12 & 
   4.8e-06 ,3.45e-07
  \\
  SubApSnap-LU & 0.302 & 1.6e-07 &      1.88e-05 ,  1.7397e-06      \\
  SubApSnap-Lev & 0.469 &  3.6e-11   &  5.8e-06  ,  4.12e-07   \quad 
  \end{tabular}
\end{table}

\subsection{Kernel ridge regression with hyperparameter optimisation}

In this example we consider using SubApSnap for 
hyper-parameter tuning in the context of 
Kernel Ridge Regression (KRR) applied to noisy data. 
The goal is to find a smooth function that approximates the original signal from noisy training samples using an RBF kernel. 

Specifically,
we generate $\hat n = 11,000$ equispaced (but randomly permuted) points $\{t_1, \ldots, t_{\hat n}\}$ in the interval $[0,10]$. 
The targets are defined as
\begin{equation}
y_i = \sin(t_i) + \epsilon_i, \qquad \epsilon_i \sim \mathcal{N}(0, 0.3^2).
\end{equation}
That is, the original function is $y=\sin(t)$, and the data comes with iid noise. 
We split the first $n=10,000$ samples for training and use the last $1000$ for testing.

Given parameter values $\sigma$ (length scale of the kernel) and $\lambda$ (regularisation parameter), 
along with 
training data $(t_{\text{train}}, y_{\text{train}})$, KRR solves the $n\times n$ linear system 
\begin{equation}\label{eq:KRReqn}
(K + \lambda I) x =  y_{\text{train}},
\end{equation}
where $K \in \mathbb{R}^{n \times n}$ is the kernel matrix with entries
\begin{equation}\label{eq:K}
K_{ij} = k(t_i, t_j) = \exp\!\left(-\frac{\|t_i - t_j\|_2^2}{2\sigma^2}\right), 
\end{equation}
and $\lambda > 0$ is the ridge regularization parameter.

Once the vector $\alpha\in\mathbb{R}^n$ is obtained, the function value at a 
new point $t$ is approximated by 
$f(x) = k(t, T_{\text{train}})^T x$, where $T_{\text{train}}=[t_1,\ldots,t_n]^T$.
Clearly, the outcome depends on the hyperparameters $(\lambda,\sigma)$, for which the best values are unknown. 
Here we tune the hyperparameters $\lambda$ and $\sigma$ by grid search. 
That is, we take 30 logarithmically spaced points $\lambda_i$ where $\lambda_1=10^{-5},\lambda_{20}=10^2$, and 30 equispaced points $0.1=\sigma_1,\sigma_2,\ldots,\sigma_{20}=10$, and form a $20\times 20$ grid of $(\lambda_i,\sigma_j)$ values. 
For each of these 900 pairs, we train on the training set and compute the root-mean-squared error (RMSE) on the test set. The pair $(\lambda, \sigma)$ achieving the lowest RMSE is selected as the hyperparameter. We refer to~\cite{rasmussen2006gaussian,scholkopf2002learning,shawe2004kernel} for details on kernel methods. 

The main computational task in this process is in 
(i) forming the kernel matrix $K$ (which clearly depends on $\sigma$) using~\eqref{eq:K}, and (ii) solving the linear system~\eqref{eq:KRReqn}. 
Note that SubApSnap is able to reduce the cost dramatically for both. That is, once the snapshot and subsampling matrices are chosen, (i) reduces to evaluating $K$ at only the row-subset, and (ii) reduces to computing $SB$ and solving an $O(r)\times r$ least-squares problem. 
The reduction in cost for (i) can be as large as a factor $n/r$, and (ii) is as described in Tables~\ref{tab:comlexity},\ref{tab:complexity_newp}. 

\begin{figure}
  \begin{minipage}[t]{0.5\hsize}
          \includegraphics[width=0.93\linewidth]{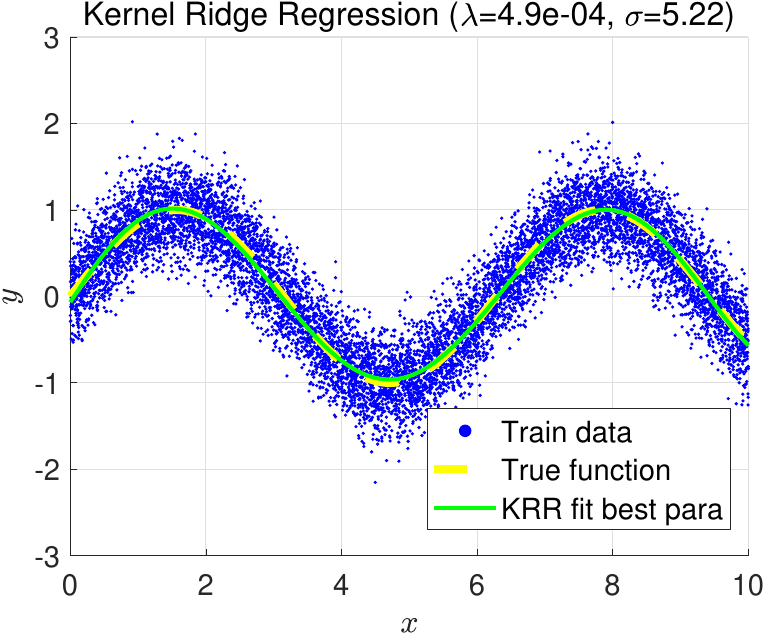}
  \end{minipage}   
  \begin{minipage}[t]{0.5\hsize}
      \includegraphics[width=0.97\linewidth]{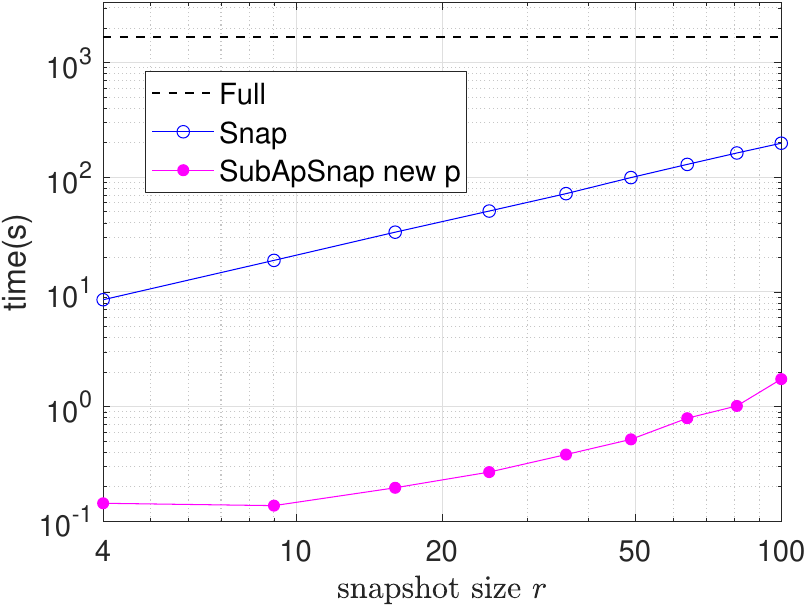}  
  \end{minipage}  
    \caption{
    Kernel Ridge regression example. 
        Left: original function, data points, and the approximation computed by KRR with the tuned hyperparameters $(\lambda,\sigma)$. With a randomly chosen $(\lambda,\sigma)$, the error can be significantly larger, e.g. $>1$. 
    Right: Runtime of SubApSnap (where the time for computing the snapshots ($O(n^3r)$) is separately reported from 
    generating and solving    
    $\min_{\hat c}\|S(A(p_*)X\hat c-b)\|_2$, which is $O(nr^2\ell)=O(900nr)$).
    }
    \label{fig:KRRtime}
\end{figure}

The results are reported in Figure~\ref{fig:KRRtime}. 
In this example, the snapshot rank $r$ is a squared number, as we take $\sqrt{r}\times \sqrt{r}$ grid points in the $(\lambda,\sigma)$ space.
We saw that $r> 50$ was required to find the optimal parameter values. 
We see that with such $r$, for each new value of $p$, SubApSnap is achieving speedup of more than three orders of magnitude. 
In this example, SubApSnap's relative residual was bounded by $10^{-8}$ for all $r\geq 16$, and their geometric mean was $5.4\times 10^{-12}$. 

\section{Discussions and other possible applications}
SubApSnap combines two simple ideas to provide a powerful approach to solving many related linear systems. In our medium-to-large-scaled computations it is already producing significant speedup; in extreme-scale problems, the benefit is certain to be more dramatic. 

We close with a few more examples where SubApSnap can be helpful. 

\paragraph{Least-squares problems.} 
The key ideas of SubApSnap are directly applicable to parameter-dependent least-squares problems $\min_x\|A(p)x-b(p)\|_2$: compute a snapshot consisting of minimisers $x_i$ of $\|A(p_i)x-b(p_i)\|_2$, form a subsampling matrix $S$, and solve $\min_{\hat c}\|S(A(p_*)X\hat c-b(p_*))\|_2$. 

\paragraph{Nonlinear eigenvalue problems.} 
Given a matrix-valued function $N(\lambda)\in\mathbb{C}^{n\times n}[\lambda]$, 
the goal of nonlinear eigenvalue problems is to find $\lambda\in\mathbb{C}$ and $x\in\mathbb{C}^n$
 $N(\lambda)x=0$. One approach is to find the poles of the resolvent $H(s) = u^T(N(s)^{-1})v$~\cite{brennan2023contour,bruno2024evaluation}. To do so we need to sample $H$ at many values of $s$; each of these requires a linear system solution $N(s_i)x_i=v$. Clearly, one can use SubApSnap for this. 

\paragraph{Newton's and gradient methods.} 
As mentioned in the introduction, 
Newton's method for minimising $f(x)$ reduces to  $A(p)$ is the Hessian of $f$ at $p=x\in\mathbb{R}^n$. 
Similarly, Newton's method for zerofinding $f(x)=0$ where $f:\mathbb{R}^n\rightarrow \mathbb{R}^n$ requires a linear system wherein the matrix is the Jacobian $\nabla f(x)$, which usually depends smoothly on the input $x$.  
Alternatively, in a gradient-based method, one could form a snapshot consisting of the gradient vectors at previous iterations. 

\paragraph{Parameter-dependent optimisation.} 
Another situation in optimisation is when minimisation is required for a parameter-dependent function $f(x,p)$ for many parameter values $p$~\cite{bonnans2013perturbation}. In lieu of solving the Newton equation, one can find a SubApSnap solution using a snapshot based on solutions for minimising the "snapshot functions" $f(x,p_i)$. 

In view of the steady growth of applications in scientific computing and data science that reduce to linear systems, it is safe to assume that parameter-dependent linear systems will continue to grow in number and importance, where SubApSnap may be able to dramatically speed up the solution.

\bibliographystyle{abbrv}
\bibliography{bib2}

\end{document}